\documentclass[10pt,a4paper]{article}

\usepackage{amssymb,amsbsy,amsmath,amsfonts,amssymb,amscd, mathrsfs}
\usepackage{enumerate}

\usepackage{color}
\usepackage[english]{babel}
\usepackage[T1]{fontenc}
\usepackage{indentfirst}
\makeatletter
\@addtoreset{equation}{section}
\makeatother

\newtheorem{statement}{}[section]
\newtheorem{definition}[statement]{Definition}
\newtheorem{theorem}[statement]{Theorem}
\newtheorem{proposition}[statement]{Proposition}
\newtheorem{corol}[statement]{Corollary}
\newtheorem{lemma}[statement]{Lemma}

\newtheorem{ex}[statement]{Example}

\newcommand{\C}{\mathbb C}
\newcommand{\Prob}{\mathbb P}
\newcommand{\D}{\mathbb D}
\newcommand{\R}{\mathbb R}

\newcommand{\N}{\mathbb N}
\newcommand{\T}{\mathbb T}
\newcommand{\E}{\mathbb E}
\newcommand\e{{\rm e}}
\newcommand{\eps}{\varepsilon}
\newcommand\ind{{\rm 1\kern-.30em I}}

\newcommand\dis{\displaystyle}
\newcommand{\cqfd}
{%
\mbox{}%
\nolinebreak%
\hskip7pt
\rule{2mm}{2mm}%
\medbreak%
\par%
}

\newcommand{\biindice}[3]%
{%

\begin{array}[t]{c}
{\displaystyle #1}\\
{\scriptstyle #2}\\
{\scriptstyle #3}
\end{array}

}

\begin{document}

\title{\bf Absolutely summing Carleson embeddings on Hardy spaces}
\author{Pascal Lef\`evre and Luis Rodr\'\i guez-Piazza\footnote{Partially supported by the project MTM2015-63699-P (Spanish MINECO and FEDER  funds)}}
\date{ \footnotesize \today}

\maketitle

\bigskip

 {\bf Abstract.} {\it We consider the Carleson embeddings of the classical Hardy spaces (on the disk) 
 into a $L^p(\mu)$ space, where $\mu$ is a Carleson measure on the unit disk. 
 This includes the case of composition operators. We characterize such operators 
 which are $r$-summing on $H^p$, where $p>1$ and $r\ge1$. This completely extends the former results on the subject and solves a problem open since the early seventies.}
\medskip

\par\medskip


\section{Introduction} 

In this paper, we investigate Carleson embeddings on classical Hardy spaces 
$H^p$ when $p>1$. In the following, the unit disk of the complex plane is 
denoted $\dis\D=\big\{z\in \C\big|\, |z|<1\big\}$. Its boundary, the torus, is 
denoted $\dis\T=\big\{z\in \C\big|\, |z|=1\big\}=\partial\D$. 
We denote by ${\cal H}(\D)$ the class of holomorphic functions on the unit disk. 
At last the Hardy spaces are defined by
$$\dis H^p=\bigg\{f\in {\cal H}(\D)\big|\,\sup_{r<1}\int_{\mathbb T}\big|f\big(rz\big)\big|^p\,d\lambda<\infty\bigg\}$$

and
  $${\|f\|_p=\sup_{r<1}\Big(\int_{\mathbb T}|f(rz)|^p\,d\lambda\Big)^{1/p}} =\sup_{r<1}\|f_r\|_{L^p({\mathbb T})}\;.$$
Here $\lambda$ stands for the normalized Haar measure on the torus 
(it is the normalized arc length), and $f_r(z)=f(rz)$ with $r\in(0,1)$ and $z\in\overline{\D}$.
\medskip

Now, let us turn to our main subject. Given a positive Borel measure $\mu$ on the closed 
unit disk $\overline{\mathbb D}$, we consider the formal identity $J_\mu$ 
from the Hardy space $H^p$ into $L^p(\mu)$ (we keep the notation $J_\mu$ 
instead of $J_{p,\mu}$ in the sequel for sake of lightness):

\medskip

\centerline{
\begin{tabular}{cccl}
$J_\mu$:& $H^p$& $\longrightarrow $& $L^p(\mu)$\cr
& $f$ & $\longmapsto$ & $f$\cr
\end{tabular}}
\medskip

Thanks to a famous result of Carleson (see \cite{C}), this is well defined and bounded if and only if $\mu$ is a Carleson measure, i.e. 
$$\displaystyle\sup_{\xi\in\mathbb T}\mu\Big({\cal W}(\xi,h)\Big)=O(h),\qquad\hbox{when }h\to 0,$$
where ${\cal W}(\xi,h)$ is the Carleson window 
$${\cal W}(\xi,h)=\{z\in{\mathbb C}\;|\, 1-h\le|z|\le 1\quad\hbox{and}\quad|\arg(z\bar\xi)|\le h\}.$$

Let us recall that $J_\mu$ is compact if and only if $\mu$ is a vanishing Carleson measure (see \cite{Po}, or \cite{McC} for composition operators in higher dimension): 
$$\displaystyle\sup_{\xi\in\mathbb T}\mu\Big({\cal W}(\xi,h)\Big)=o(h),\qquad\hbox{when }h\to 0.$$

Moreover let us mention that there is no real restriction in assuming in the sequel 
of this paper that $\mu$ is actually a measure carried by the open unit disk: let us 
assume {\sl a priori} that $\mu$ is carried by the closed unit disk and that $J_\mu$ is either 
order bounded or $r$-summing for some $r\ge1$, then in both cases, $J_\mu$ is actually 
compact and a necessary condition is that $\mu({\mathbb T})=0$. Hence from now 
on to the end of the paper, we are going to assume that $\mu$ is actually a positive 
Borel measure $\mu$ on the open unit disk ${\mathbb D}$.
\medskip

One motivation to get interested in the Carleson embedding is that it allows 
to treat the case of composition operators on $H^p$. 
\smallskip

Let us recall that, given  a symbol, {\sl i.e.} an analytic function 
$\varphi:\D\to\D$, the composition operator $C_\varphi:H^p\to H^p$ is well defined and 
automatically bounded (see the monographs \cite{CmC} or \cite{S1} for example). 
Moreover many operator properties of $C_\varphi$ can be expressed in terms of Carleson measures 
thanks to the transfert formula. Indeed, since $\varphi\in H^\infty$, 
it admits boundary values almost everywhere: $\dis\lim_{r\rightarrow1^-}\varphi(r\xi)$ 
exists for almost every $\xi\in\T$. We shall write simply $\varphi(\xi)$ in the sequel 
(although in the literature, it is often denoted $\varphi^\ast(\xi)$). 
The  pullback measure of $\lambda$ associated to $\varphi$ plays now a crucial role:
$$\lambda_\varphi(E)=\lambda\Big(\big\{\xi\in\T|\; \varphi(\xi)\in E\big\}\Big)\quad\hbox{for every  Borel subsets } 
E\hbox{ of }\overline{\D}.$$

The transfer formula gives

$$\|f\circ\varphi\|_{H^p}=\|f\|_{L^p(\D,\lambda_\varphi)}\qquad\hbox{for every } f\in H^p\,.$$
Hence many properties of the operator $C_\varphi$ are common with 
the ones of the operator $J_{\lambda_\varphi}$, in particular compactness, $r$-summingness,... 

\smallskip

The case of weighted composition operators can also be treated in the same manner.
\medskip

We are going to characterize those operators $J_\mu$, 
which are $r$-summing for some $r\ge1$. Before standing the results, let us recall the definitions

\begin{definition}
Suppose $1\le r<+\infty$ and let $T\colon X\to Y$ be a (bounded) operator between 
Banach spaces. We say that $T$ is an  {$r$-summing operator}  if  there exists $C\ge0$ such that 
$$
\Big(\sum_{j=1}^n \|Tx_j\|^r \Big)^{1/r}\le 
C \sup_{x^*\in B_{X^*}}\Big(\sum_{j=1}^n |\langle x^*, x_j\rangle|^r\Big)^{1/r}=  
C \sup_{a\in B_{\ell^{r'}}}\Big\|\sum_{j=1}^n a_jx_j\Big\|
 $$
for every finite sequence $x_1$, $x_2,\ldots,x_n$ in $X$. 
\smallskip

The $r$-summing norm of $T$, denoted by $\pi_r(T)$, is the least suitable constant $C\ge0$. 
\end{definition}

The class of $r$-summing operators forms an operator ideal (for instance see \cite{DJT} for more details). 
\medskip

We shall use several times the following well known fact about summing operators (see \cite{DJT}): 
a bounded operator $T:X\rightarrow Y$ is $r$-summing if and only if there exists $C>0$ such that, 
for every $X$-valued random variable $F$ on any measure space $(\Omega,\nu)$, we have 

\begin{equation}\label{rsummingintegral}
\int_\Omega\|T\circ F\|^r\,d\nu\le C\sup_{\xi\in B_{X^\ast}}\int_\Omega|\xi\circ  F|^r\,d\nu
\end{equation}
Actually the best admissible $C$ is $(\pi_r(T))^r$.

\medskip

Very few results are known on absolutely summing  composition operators: 
there is a characterization of $r$-summing composition operators on $H^p$ due to 
Shapiro and Taylor in \cite{ST} only  when $r=p\ge2$. 
The same result (with an obviously adapted proof) is actually valid for general Carleson embeddings:

\begin{theorem}{\cite{ST}}
Let $p\ge2$.

$J_\mu$ is $p$-summing on $H^p$ if and only if

\begin{equation} \label{ordreborne}
\dis\int_{\mathbb D}\frac{1}{1-|z|}\, d\mu <\infty. \qquad\qquad
\end{equation} 
\end{theorem}

Moreover, for every $p\ge 1$, the condition (\ref{ordreborne}) is sufficient 
to ensure that $J_\mu$ is a $p\,$-summing operator on $H^p$. When $1\le p\le2$, $J_\mu$ is 
actually even absolutely summing since $H^p$ has cotype $2$.

A natural question then arises: is (\ref{ordreborne}) the good condition (i.e. a necessary condition) 
when $1\le p\le2$? This is false in general: Domenig proved in \cite{Do} that, given $p\in[1,2)$ 
there exists an absolutely summing composition operator on $H^p$ which is not order bounded. 
He was able to give some sufficient condition for the construction of his example, 
but without any characterization. Let us mention that, in this case, 
it is equivalent to be an order bounded Carleson embedding and to verify 
condition (\ref{ordreborne}). Indeed, the following is  known from the 
specialists (see for instance \cite{LLQR1} in the case of composition operators).

\begin{proposition}
Let $\mu$ be a Carleson measure on the open unit disk $\D$ and $p\ge1$.

$J_\mu\colon H^p\rightarrow L^p(\mu)$ is order bounded if and only if (\ref{ordreborne}) is satisfied.
\end{proposition}

\noindent{\bf Proof.} By definition, $J_\mu$ is order bounded if and only if 
there exists some $h\in L^p(\mu)$ such that for every $f$ in the unit ball of $H^p$, we have 
$\dis |f|\le h$ a.e. on $\D$. Since $H^p$ is separable, it suffices to test this control on a 
dense countable subset of the unit ball de $H^p$. Hence $J_\mu$ is order bounded if and only if
$$\dis\int_{\mathbb D}\biindice{\sup}{f\in H^p}{\|f\|\le1}|f(z)|^p\;d\mu<\infty$$
which is equivalent to
$$\dis\int_{\mathbb D}\|\delta_z\|^p\;d\mu<\infty,$$
where $\delta_z$ is the evaluation at the point $z\in{\mathbb D}$, viewed as a functional on $H^p$. 

It is well known that $\|\delta_z\|=\dis\frac{1}{(1-|z|^2)^{1/p}}$ and the result follows.\cqfd

\medskip

As far as we know, there is no other result on the characterization of Carleson embeddings (or merely summing composition operators). In particular, the following problem was fully open 
(except when $r=p\ge2$).
\medskip

{\bf Problem:}

{\sl Given $p,r\ge1$ and $\mu$ a Carleson measure on the unit disk, 
which condition on $\mu$ characterizes the fact that  $J_\mu$ is a $r$-summing operator~?}
\medskip

In this paper, we are going to solve completely this problem for $p>1$ and $r\ge1$. 
Our characterizations involve Carleson windows when $p\ge2$ and $r>p'$, 
or integral conditions  when $p\in(1,2)$ or $r\le p'$. Here, as usual, $p'$ is 
the conjugate exponent of $p\ge1$: $\dis\frac{1}{p}+ \frac{1}{p'}=1.$

Concretely, let us describe the organization of the paper. This first section is devoted to the introduction of the notions and questions. In the second one, we state our main results and specify the special case of composition operators. We finish with a few examples enlightening our statements. The third one treats what we call the diagonal case: we restrict the domain to the space spanned by the monomials $z^{N+1},\ldots,z^{2N}$ and the measure to the corona $|z|\in[1-1/N,1-1/2N)$. It turns out that in this framework, the embedding acts as a diagonal operator on the classical $\ell^p$ space. The difficulty is then to glue the pieces: section 4 provides some tools to do so and section 5 makes explicits some consequences for our purpose. In section 6, we focus on the case $p>2$: most of the cases follow the results obtained in the previous sections, but the case $r\le p'$ requires a specific approach. The case $p<2$ cannot be treated following the same ideas and we treat it in the section 7. At last section 8 is devoted to some examples and remarks. For instance, we focus on the (formal) identity from the Hardy space $H^p$ to the Bergman space ${\cal B}^q$, and on the other hand we compare the different classes of $r$-summing  composition operators.

Results of this paper were announced (without proof) in \cite{LR}.

\medskip

We state here a first statement which is an elementary necessary condition when $1\le p\le2$, 
and has a double interest: it first gives a useful (practical) test in some cases. 
On the other hand, this necessary condition is a step, which turns out to be mandatory 
for a step in the proof of our characterization. Without waiting for the right characterization, 
we can notice that this is surely not a general sufficient condition.

\begin{proposition}\label{propnec} Let $1\le p\le 2$.

\begin{enumerate}
\item Let $1\le q\le 2$. We assume that the formal identity 
$f\in H^p\longmapsto L^q(\D,\mu)$ is a $r$-summing operator for some $r\ge1$. 
Then the measure $\dis\frac{d\mu}{(1-|z|)^\frac{q}{2}}$ is finite:
\begin{equation} \label{ncpq}
\dis\int_{\mathbb{D}}\frac{1}{(1-|z|)^\frac{q}{2}}\; d\mu<\infty.
\end{equation} 

In particular, if we assume that $J_\mu:H^p\to L^p(\mu)$ is an $r$-summing operator for some $r$ then
\begin{equation} \label{ncp}
\dis\int_{\mathbb{D}}\frac{1}{(1-|z|)^{p/2}}\; d\mu<\infty.
\end{equation} 

\item When $1\le p< 2$, the preceding condition~\ref{ncp} is not sufficient in general.
\end{enumerate}
\end{proposition}

For instance, applying this result to the normalized (area) measure ${\cal A}$ on $\D$, 
we can already point out that the formal identity from $H^p$ to the Hilbert Bergman space 
is not $r$-summing for any $p\in[1,2]$ and any $r\ge1$.

{\bf Proof.} 1. Since $H^p$ and $L^q(\mu)$ have cotype $2$, the operator is actually 
$r$-summing for every $r\ge1$. In particular, it is a $q$-summing operator. 
Thanks to the Pietsch domination theorem, there exist $C>0$ and a probability 
measure $\nu$ on the unit ball of $(H^p)^\ast$ such that, for every $f\in H^p$, we have

$$\dis\int_{\mathbb{D}}|f|^q\, d\mu\le C\displaystyle\int_{B_{(H^p)^\ast}}|\alpha(f)|^q\, d\nu(\alpha).$$

We apply this inequality to $$f_\omega(z)=\dis\sum_{n=0}^Nr_n(\omega)z^n$$ 
where $N\ge1$ and $(r_n)$ is a sequence of i.i.d. random variables  (Rademacher). 
Then we can take the expectation with respect to $\omega$ to obtain via Fubini and the Khinchin inequalities:  

$$\dis\int_{\mathbb{D}}\Big(\dis\sum_{n=0}^N|z|^{2n}\Big)^{q/2}\, d\mu\le 
C'\int_{B_{(H^p)^\ast}}\E|\alpha(f_\omega)|^q\, d\nu(\alpha).$$
for some constant $C'$ depending only on $q$ and $C$.

But, for any $\dis\alpha\in B_{(H^p)^\ast}$, we can write $\alpha(f)=\int_{\mathbb T}g\big(\overline{z}\big) \,f(z) d\lambda$, where $g$ belongs to the unit ball of $L^{p'}$, so that (recall that $q\le2$)
$$\E|\alpha(f_\omega)|^q\le\Big(\E|\alpha(f_\omega)|^2\Big)^{q/2}=\Big(\E\big|\sum_{n=0}^Nr_n(\omega)\widehat{g}(n)\big|^2\Big)^{q/2}\le\|g\|_2^q\le1\; ,$$
as $p'\ge2$ and every $g\in B_{L^{p'}}$ belongs to $B_{L^2}$. 

We get for arbitrary large $N$:
$$\dis\int_{\mathbb{D}}\Big(\frac{1-|z|^{2(N+1)}}{1-|z|^2}\Big)^\frac{q}{2}\, d\mu\le C'.$$ 
Taking the limit when $N\rightarrow+\infty$, the conclusion follows.

2. We consider the examples contructed in Lemma 4.3.\cite{LLQR3}. 
There exists some function $\varphi:\D\to\D$ analytic where 
$\big|\varphi\big(\e^{it}\big)\big|={\rm e}^{- f(t)}$ with $f(t)\sim|t|$ for $t$ 
in the neighborhood of $0$ and $\big|\varphi\big(\e^{it}\big)\big|<1$ out of this neighborhood of $0$. 
Since $p/2<1$, it is clear that condition~\ref{ncp} is fulfilled although $C_\varphi$ is not summing on $H^p$: 
indeed it would be compact on $H^p$ (thanks to Sarason's result in \cite{Sa} for $p=1$), equivalently, 
compact on $H^2$ which is not.\cqfd
\medskip
\section{The main results.}

We state below our main results which characterize any absolutely summing Carleson embedding. 
As usual, the notation $A\approx B$ means that there exist two constants $c,c'>0$
 (depending on $r$ and $p$ only) such that $A\le cB\le c'A$.

In addition to the Carleson windows, these characterizations involve special domains. 
We first divide the open unit disk $\D$ into dyadic annuli 
$$
\Gamma_n=\big\{z\in\D\; \big| \, 1-\frac{1}{2^n}\le|z|<1-\frac{1}{2^{n+1}}\big\},\quad\hbox{where }n=0,1,2,\cdots
$$
Then each $\Gamma_n$ is divided into $2^n$ similar pieces $R_{n,j}$, $0\le j<2^n$,
that we call Luecking boxes (or Luecking rectangles): 
$$
R_{n,j}=\big\{z\in\D\,  \big|\; 1-\frac{1}{2^n}\le|z|<1-\frac{1}{2^{n+1}}
\hbox{ and }\arg(z)\in\big(2\pi j/2^n,2\pi(j+1)/2^n  \big]\big\}.
$$
So clearly the family of the $R_{n,j}$ with $n\ge0$ and  $0\le j<2^n$, forms a partition of $\D$.
\medskip

We will also use the Stolz domain $\Sigma_\xi$ at $\xi\in\T$ which 
is the interior of the convex hull of $D(0,1/2)\cup\{\xi\}$.

\bigskip

{\bf Characterization of absolutely summing Carleson embeddings:}
\medskip

{\sl Let $\mu$ be a Carleson measure on $\D$. 
\medskip

1) Let $1<p\le2$. The natural injection 
$J_\mu\colon H^p(\D)\to L^p(\mu)$ is $2$-summing if and only if 
$$
\Phi\colon 
\xi\longmapsto\int_{\Sigma_\xi}  \frac{1}{\big(1-|z|^2\big)^{1+\frac{p}{ 2}}}d\mu(z)
$$ 
belongs to $L^{2/p}(\T,d\lambda)$, where $\Sigma_\xi$ is the Stolz domain at point $\xi\in\T$. Moreover we have
\begin{equation}\label{p<2}
\pi_2(J_\mu)\approx \|\Phi\|_{2/p}^{1/p}\approx \|\Psi\|_{2/p}^{1/p}\,,
\end{equation}where $\dis\Psi(\xi)=\int_\D\frac{1}{|1-\overline{\xi}z|^{1+\frac{p}{ 2}}}\,d\mu(z)$, for every $\xi\in \T$.

\bigskip

2) Let $p\ge 2$ and $r\ge 1$.
\medskip

\begin{itemize}
\item When $1\le r\le p'$ , we have
\begin{equation}  \label{1<r<p'}\pi_r(J_\mu) \approx\pi_1(J_\mu) \approx
\dis \|F\|_{L^{p'}({\mathbb T})},\quad\hbox{where }F(\xi)=\dis\Bigg[\sum_{n\ge0}2^{2n}
\Big(\mu\big(W(\xi,2^{-n}\big)\Big)^\frac{2}{p}\Bigg]^{1/2}\end{equation}

\item When $p'< r\le p$, we have
\begin{equation}\label{p'<r<p}
\pi_r(J_\mu)\approx
\dis\Bigg[\sum_{n\ge0}\;\sum_{0\le j<2^n}\Big(2^{n}\mu(R_{n,j})\Big)^{r/p}\Bigg]^{1/r}
\end{equation}

\item  When $p\le r$, we have
\begin{equation}\label{r>p}\pi_r(J_\mu) \approx\pi_p(J_\mu)\approx
\Bigg[\sum_{n\ge0}\;\sum_{0\le j<2^n} 2^{n}\mu(R_{n,j})\Bigg]^{1/p}\approx
\Bigg(\int_{\D}\frac{1}{1-|z|}\;d\mu\Bigg)^{1/p}\end{equation}
\end{itemize}
}
\medskip

Before giving the corollaries for composition operators, let us mention several remarks:

\begin{enumerate}[i)]
\item In the case $p=2$, we recover the characterization of $2$-summing operators (i.e. Hilbert-Schmidt operators): 
$$
\int_{\mathbb T}\int_{\Sigma_\xi}\dis\frac{1}{(1-|z|)^2}d\mu\; d\lambda(\xi)=
\int_{\mathbb D}\int_{\{\xi|\ z\in\Sigma_\xi\}}\dis\frac{1}{(1-|z|)^2}d\lambda(\xi)\;d\mu\approx
\int_{\mathbb D}\dis\frac{1}{1-|z|}\;d\mu  $$
thanks to Fubini's theorem and since $\dis\{\xi|\ z\in\Sigma_\xi\}$ is an interval of length $\approx 1-|z|$.
\medskip

\item On the other hand, we recover the necessary condition given by Prop~\ref{propnec}: when $p\le2$, we have $2/p\ge 1$, and hence, by Fubini's theorem, 
$$
\begin{array}{ccl}
\dis
\Bigg(\int_{\mathbb T}\Bigg(\int_{\Sigma_\xi}\dis\frac{1}{(1-|z|^2)^{1+\frac{p}{2}}}
d\mu\;\Bigg)^{2/p} d\lambda \Bigg)^{p/2}&\ge&\dis
\int_{\mathbb T}\int_{\Sigma_\xi}\dis\frac{1}{(1-|z|^2)^{1+\frac{p}{2}}} \;d\mu d\lambda\cr
&=&\dis\int_{\mathbb D}\int_{\{\xi|\, z\in\Sigma_\xi\}}\dis\frac{d\lambda}{(1-|z|^2)^{1+\frac{p}{2}}} \;d\mu\cr 
&\approx&\dis\int_{\mathbb D}\dis\frac{1}{(1-|z|^2)^\frac{p}{2}}d\mu\cr
\end{array}$$
\medskip

\item Our characterizations show that the $r$-summing character of $J_\mu$ depends only on the sequence of values $\big\{\mu\big(R_{n,j}\big)\big\}_{n\ge0;j<2^n}$. More precisely, we point out that when two positive (finite) measures $\mu$ and $\nu$  satisfy $\dis\mu\big(R_{n,j}\big)\le\nu\big(R_{n,j}\big)$ for every Luecking rectangle $R_{n,j}$, then $\dis\pi_r\big(J_\mu\big)\lesssim\pi_r\big(J_\nu\big)$ for any $r\ge1$. This can be checked from our characterizations: it is obvious when $p\ge2$ and $r>p'$. In the other cases, just use the fact that the Luecking rectangles forms a partition of the unit disc. For instance, the function $\Psi$ in (\ref{p<2}) is equivalent to
$$\biindice{\sum}{n\ge0}{j<2^n}d\big(\xi,R_{n,j}\big)^{-(1+p/2)}\mu\big(R_{n,j}\big).$$

In particular, when $\dis\mu\big(R_{n,j}\big)=\nu\big(R_{n,j}\big)$ for every $n$ and $j$, then $J_\mu$ is $r$-summing if and only if $J_\nu$ is $r$-summing.

\medskip

\item This leads to the natural question to wonder whether the characterizations depend on the order of enumeration of the values $\mu\big(R_{n,j}\big)$. More precisely, when $p\ge2$ and $r>p'$, the $r$-summing character of $J_\mu$ is clearly invariant by permutation of the values of $\dis\big\{\mu\big(R_{n,j}\big)\big\}_{j<2^n}$ (for each fixed integer $n$). It turns out that it is no more true in the other cases. We have the following examples
\begin{itemize} 
\item {\bf Example 1.} Let $p>2$. There exist two (finite) measures $\mu$ and $\nu$ on $\D$ such that 
\begin{itemize}

\item For every $n\ge1$, the sequence $\dis\big\{\mu\big(R_{n,j}\big)\big\}_{j<2^n}$ is a permutation of the sequence $\dis\big\{\nu\big(R_{n,j}\big)\big\}_{j<2^n}$.

\item $J_\mu:H^p\rightarrow L^p(\mu)$ is $1$-summing.

\item $J_\nu:H^p\rightarrow L^p(\nu)$ is not $p'$-summing.
\end{itemize}
\smallskip

Consider the centers $z_{n,j}$ of the $R_{n,j}$ and two sequences: $(m_n)$ is defined as the integer part of $\dis\frac{2^n}{n^2}$ and $\alpha_n=\dis 2^{-np}n^{\gamma p}$ with 
$\gamma$ fixed in the interval $\big(\frac{1}{p'},\frac{2}{p'}-\frac{1}{2}\big)$.

Now our measures are
$$\mu=\sum_{n>n_0}\alpha_n\mu_n\quad\hbox{where}\quad \mu_n=\alpha_n\sum_{j=0}^{m_n}\delta_{z_{n,j}}$$
and 
$$\nu=\sum_{n>n_0}\alpha_n\nu_n\quad\hbox{where}\quad \nu_n=\alpha_n\sum_{j=l_n}^{l_n+m_n}\delta_{z_{n,j}}$$
where $n_0$ is large enough and $\dis l_{n+1}=2\big(l_n+m_n\big)$.

Let us point out that the radial projection on $\T$ of the rectangles $R_{n,j}$ charged by $\nu$ are pairwise disjoints arcs whereas for the measure $\mu$, these projections are arcs tending to the point $1$, when $n$ tends to infinity.


\smallskip

\item {\bf Example 2.} Let $p\in(1,2)$. There exist two (finite) measures $\mu$ and $\nu$ on $\D$ such that 
\begin{itemize}

\item For every $n\ge1$, the sequence $\dis\big\{\mu\big(R_{n,j}\big)\big\}_{j<2^n}$ is a permutation of the sequence $\dis\big\{\nu\big(R_{n,j}\big)\big\}_{j<2^n}$.

\item $J_\mu:H^p\rightarrow L^p(\mu)$ is $1$-summing.

\item $J_\nu:H^p\rightarrow L^p(\nu)$ is not $r$-summing for any $r\ge1$.
\end{itemize}

The construction is similar to the previous one but we have to interchange the way we define $\mu$ and $\nu$ in order to get now that the radial projection on $\T$ of the rectangles $R_{n,j}$ charged by $\mu$ are pairwise disjoints arcs whereas for the measure $\nu$, these projections are arcs tending to the point $1$. Also the parameters $\alpha_n$ have to be adapted.

\end{itemize}
We leave the check of the details to the reader.

\end{enumerate}
\bigskip

As an immediate corollary, we can characterize $r$-summing weighted composition 
operator for any $r\ge1$ and any $p>1$. Indeed, dealing with the operator $f\mapsto w(f\circ\varphi)$, 
it suffices to consider the measure $\mu=\nu_\varphi$ where $d\nu=|w|^p\,d\lambda$. 
Nevertheless, we state precisely these particular results only for composition operators. 
One involves the Nevanlinna counting function:
$$
\dis N_\varphi(w)= \left\{\begin{array}{cl}\sum\limits_{\varphi(z)=w}\log \frac{1}{|z|} 
& \hbox{if}\ w\not= \varphi(0)\ \hbox{and}\ w\in \varphi({\D})\\ \\0 & \hbox{else}.\end{array}\right.
$$
where the sum runs over the roots counted with multiplicity.

\medskip

{\bf Characterization of absolutely summing composition operators:}
\medskip

{\sl Let $\varphi:{\mathbb D}\to{\mathbb D}$ be an analytic map such that its boundary values satisfy
$|\varphi(w)|<1$ for almost every $w\in\T$. If $C_\varphi$ is the composition operator from $H^p$ to $H^p$, 
we have

\medskip

1) When $1<p\le2$,  for any $r\ge1$:
$$
\pi_r(C_\varphi)\approx \pi_1(C_\varphi)\approx
\Bigg[\int_{\mathbb T}\Bigg(\int_\T
\dis\frac{1}{|1 -\overline{\xi}\varphi(w)|^{1+\frac{p}{2}}}\,d\lambda(w)\;\Bigg)^{2/p} d\lambda(\xi)\Bigg]^{1/2}
$$

2) When $p\ge2$, 
\begin{itemize}
\item for $1\le r\le p'$, 
\begin{equation}\pi_r(C_\varphi) \approx \pi_1(C_\varphi)\approx
\dis \|F\|_{L^{p'}({\mathbb T})}\quad\hbox{where }F(\xi)=
\dis\Bigg[\sum_{n\ge0}2^{2n}\Big(\lambda_\varphi\big(W(\xi,2^{-n}\big)\Big)^\frac{2}{p}\Bigg]^{1/2}\end{equation}

\item for $ p'< r\le p$, 
\begin{equation} \label{pprimer<pC}
\begin{array}{ccl}
\pi_r(C_\varphi) & \approx &\dis\Bigg[\sum_{n\ge0}\;\sum_{0\le j<2^n}
\Big(2^{n}\lambda_\varphi(R_{n,j})\Big)^{r/p}\Bigg]^{1/r} \cr

&\approx & \dis\Bigg[\int_{\D}\Bigg(\frac{N_\varphi(z)}{1-|z|^2} \Bigg)^{r/p}\frac{dA}{(1-|z|^2)^2}\Bigg]^{1/r}\cr

\end{array}\end{equation}

\item  for $p\le r$, 
\begin{equation} \label{r>pC}
\pi_r(C_\varphi)\approx\Bigg(\int_{\T}\frac{1}{1-|\varphi|}\;d\lambda\Bigg)^{1/p}\end{equation} 
\end{itemize}
}
\medskip

In the case $p'< r\le p$ (see (\ref{pprimer<pC})), the equivalence with the integral quantity comes from the fact 
that the $r$-summing norm of $C_\varphi$ (on $H^p$), thanks to Luecking's characterization \cite{LueJFA}, 
turns out to be equivalent (up to an exponent) to the $\frac{2r}{p}$-Schatten norm (on $H^2$) . 
Then we use the characterization
of the membership of $C_\varphi$ to the Schatten classes, given by
Luecking and Zhu \cite{LZ} and
involving the Nevanlinna counting function.
We could also use directly the results on equivalence between the measure 
of Carleson's windows and the Nevanlinna counting function, see \cite{LLQR2} and also \cite{EK} for a recent new proof. 


\section{The diagonal case}\label{sectiondiag}


In this section, we fix an integer $N\ge1$ and we characterize the absolutely summing norm
in the case of a measure concentrated on the $N^{th}$ corona. More precisely, we consider the restriction 
$\mu_N$ of $\mu$ to the corona ${\cal G}_N$, defined by $\mu_N(A)=\mu(A\cap{\cal G}_N)$, where
$$
{\cal G}_N=\big\{z\in\D\,\big|\, 1-\frac{1}{N}\le|z|<1-\frac{1}{2N}\big\}.
$$
We are then interested in the behavior of the operator  $J_{\mu_N}$.
It turns out that the $r$-summing norm of $J_{\mu_N}$ is equivalent to that of its restriction to
the space spanned by the monomials $z^k$, where $N\le k<2N$.

Our characterization will rely on the values of the measure of the  boxes 
$$
{\cal R}_{N,j}=\Big\{z\in\D\,\Big|\; 1-\frac{1}{N}\le|z|<1-\frac{1}{2N}\hbox{ and }
\arg(z)\in\big(2\pi j/N,2\pi(j+1)/N  \big]\Big\}\, .
$$
Of course ${\cal G}_N=\dis\cup_{0\le j <N}{\cal R}_{N,j}$. 
Later we shall use the results of this section with the dyadic version $N=2^n$ 
recovering the Luecking boxes since we have $\dis R_{n,j}={\cal R}_{2^n,j}$ and  $\Gamma_n={\cal G}_{2^n}$.

The following proposition makes the link between the behavior of the restricted Carleson 
embedding and a suitable diagonal operator on classical $\ell_N^p$ space. 
We define the multiplier operator $\dis\mathfrak{M}_\beta$ from $\ell^p_N$ to $\ell^p_N$ 
by $\dis\mathfrak{M}_\beta(\e_j)=\beta_j\e_j$ where $\beta\in\C^N$ and $\{\e_j\}_{0\le j < N}$ 
is the canonical basis of $\ell^p_N$.

\begin{proposition}\label{CourSom}
For every $r\ge1$ and $p>1$, we have
$$
\pi_r\big(J_{\mu_N}\big)\approx \pi_r\big(\mathfrak{M}_\beta\big)
$$
where $\beta_{j}=\dis\big(N \mu({\cal R}_{N,j})\big)^{1/p}$, \ for $0\le j<N$.
\end{proposition}

In the previous statement, the underlying constants depend only on $p$.
Actually we are able to prove a more general statement. Before stating it, we shall give some other definitions.
\medskip

\begin{definition}
We say that a norm $\alpha$ of operator ideal is a monotone ideal norm  if the following property is fulfilled:
whenever we have three Banach spaces $X$, $Y_1$ and $Y_2$ and two operators 
$T:X\to Y_1$ and $S:X\to Y_2$ such that $\|Tx\|\ge \|Sx\|$, for every
$x\in X$, we have $\alpha(T)\ge \alpha(S)$.
\end{definition}

All $r$-summing norms are monotone. 
If $\alpha$ is monotone, we have in particular $\alpha(T)=\alpha(j\circ T)$, for every $T:X\to Y$ 
and any isometry $j:Y\to Z$. We will need the following result about monotone norms.

\begin{lemma}\label{monotone}
Let ${\mathscr I} $ be
an operator ideal with a monotone norm $\alpha$. Let $X$, $Y$ be Banach spaces and $(\Omega,\Sigma,\mu)$
a measure space. Assume that we have defined operators $S\colon X\to Y$ and $T_\omega\colon X\to Y$, 
$\omega\in \Omega$, such  that the map $\omega\mapsto T_\omega$ is measurable 
(with values in ${\mathscr I} (X,Y)$) and we have
\begin{equation}\label{intmonotone}
\|Sx\|\le \int_\Omega \|T_\omega x\|\, d\mu(\omega)\,,\qquad\text{for every $x\in X$.}
\end{equation}
Then $\dis\alpha(S)\le \int_\Omega \alpha(T_\omega)\, d\mu(\omega)$.
\end{lemma}

{\bf Proof.} We can assume that $ \int\alpha(T_\omega)\, d\mu(\omega)<+\infty$.
Let $Z=L^1(\mu,Y)$ be the space of $Y$-valued Bochner integrable functions defined on 
$(\Omega,\Sigma,\mu)$, and define $T\colon X\to Z$ by $Tx(\omega)=T_\omega x$, $\omega\in \Omega$.
Condition \eqref{intmonotone} means that $\|Sx\|\le \|Tx\|$, for every $x\in X$. Since $\alpha$ is monotone 
we have $\alpha(S)\le \alpha(T)$. We finish because it is not difficult to see that
$$
\alpha (T)\le \int_\Omega  \alpha(T_\omega)\, d\mu(\omega)\,.
$$
In fact, this is clear if $\omega\mapsto T_\omega$ is a step function and the general case follows by density.
\hfill\cqfd

\medskip 

We denote by $\dis{H^p_N}$  the subspace of $H^p$ spanned by $\{1,\ldots,z^{N-1}\}$ 
and by $\dis\widetilde{H^p_N}=z^NH^p_N$ the space spanned by $\{z^N,\ldots,z^{2N-1}\}$. 
Clearly $\dis\widetilde{H^p_N}$ and $\dis H^p_N$ are isometrically isomorphic 
via the mapping $f\in H^p_N\mapsto z^Nf$.
\medskip

\begin{theorem}\label{CourSomGen}

Let $p>1$ and $\alpha$ be a monotone norm on an operator ideal ${\mathscr I} $. 
Let $N$ be an integer, with $N\ge1$ and consider 
 the following operators: 
\begin{enumerate}[(a)]
\item $\dis J_{\mu_N}\colon f\in H^p\longmapsto f\in L^p(\D,\mu_N)$.

\item $\dis T_{p,N}\colon f\in H^p_N\longmapsto f\in L^p(\D,\mu_N)$.

\item $\dis\widetilde{T}_{p,N}\colon f\in\widetilde{H^p_N}\longmapsto f\in L^p(\D,\mu_N)$.

\item $\dis \mathfrak{M}_\beta\colon (a_j)_{0\le j<N}\in\ell^p_N\longmapsto (\beta_{j} a_j )_{0\le j<N}\in\ell^p$, 
\quad
where $\beta_j=\big(N\mu({\cal R}_{N,j})\big)^{1/p}$.

\end{enumerate}

We have
$$\dis\alpha\big(J_{\mu_N}\big)\approx\alpha\Big(\widetilde{T}_{p,N}\Big)\approx\alpha\Big(T_{p,N}\Big)\approx \alpha\big(\mathfrak{M}_\beta\big)$$
where the underlying constants are independent of  $N$ and $\alpha$.
\end{theorem}
\medskip



We shall need some lemmas. The following one is an obvious extension of 
a classical result (see 7.10\cite{Z}, p.30). 
We shall  write $\dis\theta_j=\exp\Big(\frac{2ij\pi }{N}\Big)$ for any $j\in\{0,\ldots,N-1\}$.

\begin{lemma}\label{zyg}
Let $1<p<\infty$. Then we have, for every $f\in H^p_N$,
$$
\|f\|_{p}\approx\Bigg(\frac{1}{ N}\sum_{j=0}^{N-1}\big|f(\theta_j)\big|^p\Bigg)^{1/p}\,.
$$
\end{lemma}
Note that the estimations are uniform on $N$ and depend  only on $p$.

\begin{lemma}\label{isomellp}
Let $1<p<\infty$. 

Let us define $Q(z)=1+\cdots+z^{N-1}$ and $\dis Q_j(z)={N^{-1/p'}}Q\big(\overline{\theta_j}\, z\big)$.

Then, there exist $a_p,b_p>0$ such that for every  $u_0,\ldots,u_{N-1}\in\C$:
$$
a_p\Bigg(\sum_{j=0}^{N-1} |u_j|^p\Bigg)^{1/p}\le\Bigg\|\sum_{j=0}^{N-1} 
u_jQ_j\Bigg\|_{H^p}\le b_p\Bigg(\sum_{j=0}^{N-1} |u_j|^p\Bigg)^{1/p} 
$$

\end{lemma}

{\bf Proof.} Point out that $\dis  Q_j(\theta_l)=N^{1/p}\,\delta_{j,l}$ (the Kronecker symbol) 
for all  $j,l\in\{0,\ldots,N-1\}$. Using Lemma~\ref{zyg}, we get
$$
\dis\Bigg\|\sum_{j=0}^{N-1} u_jQ_j\Bigg\|_{H^p}^p\approx
\frac{1}{ N}\sum_{l=0}^{N-1}\Big|\sum_{j=0}^{N-1} u_jQ_j(\theta_l)\Big|^p=\sum_{l=0}^{N-1}\big|u_l\big|^p\;.
$$
\hfill\cqfd

{\bf Proof of Theorem~\ref{CourSomGen}}.

 It is obvious by restriction that $\dis\alpha\big(J_{\mu_N}\big)\ge\alpha\big( T_{p,N}\big).$
\medskip

The equivalence $\alpha\big(\widetilde{T}_{p,N}\big)\approx\alpha\big( T_{p,N}\big)$ follows from the fact that 
there is an isometric isomorphism between $\dis\widetilde{H^p_N}$ and $\dis H^p_N$  (described above) 
and the fact that for any $f\in H^p$, the $L^p(\D,\mu_N)$-norm of $f$ and $z^Nf$ are equivalent 
since $|z|^N$ belongs to $\big(\frac{1}{4},\frac{1}{\sqrt\e}\big)$ when $z\in{\cal G}_N$ (and $N>1$).
\medskip

Let us prove that $\dis\alpha\big(J_{\mu_N}\big)\lesssim\alpha\big(T_{p,N}\big)$. 
For every $f\in H^p$, we consider $y_m$ which is the projection $\pi_m(f)$ of $f$ on the space $Z_m$, 
spanned by the $z^k$ when $k$ runs over $\{mN,\ldots,(m+1)N-1\}$. 
Since $p>1$, the Riesz projection is bounded and there exists a constant $r_p>0$ depending on $p$ only 
such that $\dis\|y_m\|_{H^p}\le r_p\|f\|_{H^p}$. 
So, the operator $\pi_m\colon f\in H^p\longmapsto y_m\in Z_m\subset H^p$ is uniformly bounded by $r_p$.

Now there exists $f_m\in H^p_N$ such that $z^{mN}f_m=y_m$. Moreover the correspondence 
$y_m\in Z_m\subset H^p\leftrightarrow f_m\in  H^p_N$ defines an isometric bijection $R_m$.
We have 
$$
\dis\big\|z^{mN}f_m\big\|_{L^p(\D,\mu_N)}\le\e^{-m/2}\big\|f_m\big\|_{L^p(\D,\mu_N)}\,,
$$
for every $N\ge2$. 
It means that the operator 
$$
M_m\colon g\in L^p(\D,\mu_N)\longmapsto z^{mN}g\in L^p(\D,\mu_N)
$$
has norm less than $\dis\e^{-m/2}$.

Therefore, the sum $\sigma_\mu$ of the operators $\dis J^{(m)}=M_m\circ T_{p,N}\circ R_m\circ \pi_m$, 
acting from $H^p$ to $L^p(\D,\mu_N)$ converges. Using the ideal property, 
we have for any ideal norm $\alpha$ (and in particular for the operator norm):
$$
\dis\alpha\Big(\sum_{m\ge0} M_m\circ T_{p,N}\circ R_m\circ \pi_m\Big)\le 
\sum_{m\ge0} r_p\e^{-m/2} \alpha\big( T_{p,N}\big)\le\frac{r_p\sqrt\e}{\sqrt\e-1}\, \alpha\big( T_{p,N}\big).
$$

Since  the operators $\sigma_\mu$ and  $J_\mu$ coincide, 
for instance on polynomials (hence on $H^p$), 
the operator $J_{\mu_N}$ may therefore be written as the sum of the operators 
$J^{(m)}$, moreover 
$$
\dis\alpha\big(J_{\mu_N}\big)\le\frac{r_p\sqrt\e}{\sqrt\e-1}\, \alpha\big( T_{p,N}\big).
$$
\medskip

Let us prove that $\dis\alpha\big(\mathfrak{M}_\beta\big)\lesssim \alpha\big( T_{p,N}\big)$. 
We will assume that $\mu({\cal R}_{N,j})>0$, for every $j$, the general case can be easily deduced
from this one.
From Lemma $\ref{isomellp}$, we know that $\dis H^p_N$ is isomorphic to $\ell^p_N$ 
(with constant not depending on $N$) via the mapping
$$
\dis u\in\ell^p_N\longmapsto \Psi(u)= \sum_{j=0}^{N-1} u_jQ_j\,.
$$

On the other hand, since $\dis Q_j(z)=N^{-\frac{1}{p'}}\Big(\frac{1-z^N}{1-\overline{\theta_j}z}\Big)$, 
we have $\dis|Q_j(z)| \le N^{1/p}$, for every $z\in\D$, and 
for every $z\in {\cal R}_{N,j}$,
$$
\dis N^{1/p}\lesssim\frac{\Big( 1-\dis\frac{1}{\sqrt\e}\Big)N^{-1/p'} }{|\theta_j-z|}\le|Q_j(z)|\,.
$$
This implies that 
$$
\Big\|Q_j\ind_{{\cal R}_{N,j}}\Big\|_{L^p(\D,\mu_N)}\approx \Big(N\mu\big({\cal R}_{N,j}\big)\Big)^{1/p}.
$$

Let us consider the ``diagonal'' operator 
$$
\dis \Delta=D\circ\mathfrak{M}_\beta\circ\Psi^{-1}\colon  H^p_N\longrightarrow L^p(\D,\mu_N)
$$ 
where 
$$
D\colon u\in\ell^p_N \longmapsto  
\sum_{j=0}^{N-1}\Big(N\mu\big({\cal R}_{N,j}\big)\Big)^{-1/p}u_j Q_j\ind_{{\cal R}_{N,j}}\in L^p(\D,\mu_N)
$$

\medskip
{\bf Claim1:}\qquad $\alpha\big(\mathfrak{M}_\beta\big)\approx\alpha( \Delta)$.
\medskip

It is clear by the ideal property that 
$$
\dis\alpha( \Delta)\le\|D\|\alpha\big(\mathfrak{M}_\beta\big)\big\|\Psi^{-1}\big\|
\lesssim\alpha\big(\mathfrak{M}_\beta\big)\|D\|\,.
$$
But, for every $u\in\ell^p_N$, we have 
$$
\|D(u)\|^p_{L^p(\D,\mu_N)}=
\sum_{j=0}^{N-1}\Big(N\mu\big({\cal R}_{N,j}\big)\Big)^{-1}|u_j|^p 
\Big\|Q_j\ind_{{\cal R}_{N,j}}\Big\|^p_{L^p(\D,\mu_N)}\approx \|u\|_p^p\,,
$$
hence $\alpha( \Delta)\lesssim\alpha\big(\mathfrak{M}_\beta\big)$.

Actually since $\Psi$ is an isomorphism (with uniform constants) and $D$ 
is an isomorphism on its range, we also have 
$\alpha\big(\mathfrak{M}_\beta\big)\lesssim\alpha( \Delta)$. Claim 1 is proved.

\medskip
{\bf Claim 2:} \qquad $\alpha( \Delta) \lesssim \alpha\big( T_{p,N}\big)$.
\medskip

To prove the claim,  we shall use the underlying unconditionality, so we incorporate a random 
perturbation in some of the previous operators. More precisely, 
let us  consider a random choice of signs $\sigma=(\sigma_0,\ldots,\sigma_{N-1})\in\{\pm1\}^N$
and define 
$$
\dis \Delta_\sigma=M_\sigma\circ  T_{p,N}\circ \psi_\sigma\colon  H^p_N\longrightarrow L^p(\D,\mu_N)\,,
$$ 
where 
$$
M_\sigma\colon f\in   L^p(\D,\mu_N)\longmapsto  
\sum_{k=0}^{N-1} \Big( \sigma_k\ind_{{\cal R}_{N,k}}\Big)f\in L^p(\D,\mu_N)
$$
and $\dis \psi_\sigma$ is defined by its action on the basis $\{Q_j\}_j$:  $\dis\psi_\sigma(Q_j)= \sigma_j Q_j$.

 Clearly, for any $\sigma$, $M_\sigma$ is an isometry and $\psi_\sigma$ is an isomorphism 
 with norms not depending on $\sigma$; actually $\psi_\sigma$ is conjugated (via $\Psi$) 
 to the diagonal operator on $\ell^p_N$, associated to the $\sigma_j$, which is an isometry.

It is easy to check that for every 
$j\in\{0,\ldots,N-1\}$: 
$\dis\E_\sigma \Delta_\sigma(Q_j)= Q_j\ind_{{\cal R}_{N,j}}$. 
So we have for every $f\in   H^p_N$:
$$
\dis\E_\sigma \Delta_\sigma(f)= \Delta(f).
$$
By convexity and the properties of an ideal norm, we get that 
$\dis \alpha( \Delta)\le\E_\sigma \alpha\big(\Delta_\sigma\big)$. But 
$$
\dis\alpha\big(\Delta_\sigma\big)\le\|\Psi\|.\|\Psi^{-1}\| \alpha\big( T_{p,N}\big)\,,
$$
so
$$
\dis \alpha( \Delta)\lesssim \alpha\big( T_{p,N}\big).
$$
Claim 2 is proved and we conclude that  $\dis\alpha\big(\mathfrak{M}_\beta\big)\lesssim \alpha\big( T_{p,N}\big)$.
\medskip

At last let us prove that $\dis\alpha\big(\widetilde{T}_{p,N}\big)\lesssim\alpha\big(\mathfrak{M}_\beta\big)$. 
We again assume $\mu({\cal R}_{N,j})>0$, for every $j$.
We first concentrate our attention on the box ${\cal R}_{N,0}$ and we consider a Jordan curve 
$\gamma\subset\D$ surrounding ${\cal R}_{N,0}$ 
such that the length of $\gamma$, denoted by $\ell(\gamma)$, satisfies $\ell(\gamma)\lesssim 1/N$
and satisfying $d(z,\gamma)\gtrsim  1/N$, for every $z\in {\cal R}_{N,0}$.
By the Cauchy formula, we can write for any analytic function $f$ on $\D$ and any $z\in {\cal R}_{N,0}$:
$$
f(z)=\dis\frac{1}{2i\pi}\int_\gamma \frac{f(w)}{w-z}\;dw\,.
$$
Introducing the probability measure $d\Prob(w)=\dis\frac{1}{\ell(\gamma)}|dw|$ on $\gamma$, we obtain
$$
|f(z)|\lesssim\dis \int_\gamma |f(w)|\;d\Prob(w)\,.
$$

Now, when $z\in {\cal R}_{N,j}$, clearly $\dis|f(z)|\lesssim\dis\int_\gamma |f(\theta_jw)|d\Prob(w)$. 
Introducing the operator 
$$
U_w\colon f\in   H^p_N\longmapsto  \sum_{k=0}^{N-1} f(\theta_kw)\ind_{{\cal R}_{N,k}} \in L^p(\D,\mu_N)
$$
we have then for every $z\in{\cal G}_N$
$$
|f(z)|\lesssim\int_\gamma \big|U_w(f)(z)\big| \;d\Prob(w).
$$
In particular, 
$$
\|f\|_{ L^p(\D,\mu_N)}\lesssim \int_\gamma \big\|U_w(f) \big\|_{ L^p(\D,\mu_N)} \;d\Prob(w)
$$
and, by Lemma \ref{monotone},
$$ 
\alpha\big( T_{p,N}\big)\lesssim \int_\gamma \alpha\big(U_w \big)\;d\Prob(w).
$$

To conclude, we shall now concentrate our attention on $\dis\alpha\big(U_w \big)$, where $w\in\gamma$ is fixed. 
The operator $U_w$ is the composition of three operators: $A_w$,  $\mathfrak{M}_\beta$ and $B_w$, with 
$$
\dis A_w\colon f\in H_N^p \longmapsto \dis\big(f(\theta_j w)N^{-1/p}\big)_{0\le j<N}  \in \ell_N^p
$$ 
and 
$$
\dis B_w\colon u\in  \ell_N^p \longmapsto 
\sum_{k=0}^{N-1} u_k\big(\mu({\cal R}_{N,k})\big)^{-1/p} \ind_{{\cal R}_{N,k}} \in  L^p(\D,\mu_N).
$$

As soon as the proof that $A_w$ and $B_w$ are bounded 
(with bounds independent from $w$) will be done, 
we shall get that $\alpha\big(U_w \big)\lesssim\alpha\big(\mathfrak{M}_\beta\big)$ 
and the conclusion of our assertion.

$B_w$ is obviously isometric. Writing  $w=r\hbox{\rm e}^{ia}$ and denoting by $P_r$ the Poisson kernel and by $\tau_a$ the translation on
the circle group $\T$, we notice that 
$$
A_w(f)= \big(P_r\ast\tau_a f(\theta_j) N^{-1/p}\big)_{0\le j<N}.
$$
Since, $f\in H_N^p\mapsto P_r\ast\tau_a f\in H_N^p$ is a contraction, 
it suffices to invoke Lemma~\ref{zyg} to conclude that $A_w$ is bounded (with bound independent of $w$).
This ends the proof of the theorem.\cqfd


\section{Toolbox: how to glue summing operators}

\subsection{Summing multipliers}

The results of this section are certainly well known from the specialists. 
Nevertheless, most of them do not appear easily in the literature (actually we did not find some of them). 
For sake of completeness, we state and prove all of them.
In this section $\beta=(\beta_n)$ is a bounded sequence of complex numbers. 
In the following result, $\mathfrak{M}_\beta$ stands for the multiplier operator  on $\ell^p$, 
with $p\ge1$, defined by $\mathfrak{M}_\beta(\e_n)=\beta_n \e_n$, 
where $(\e_n)$ denotes the canonical basis of $\ell^p$.

\begin{proposition}\label{multi}
With constants only depending on $p$ and $r$, we have:

\begin{enumerate}
\item For $1\le p\le 2$ and every $r\ge1$, 
\begin{equation}\label{multi2}
\pi_r(\mathfrak{M}_\beta)\approx\|\beta\|_2
\end{equation}

\item For $p\ge 2$ and $r\le p'$, 
\begin{equation}\label{multir<p'}
\pi_r(\mathfrak{M}_\beta)\approx\|\beta\|_{p'}
\end{equation}

\item For $p\ge 2$ and $p'\le r\le p$, 
\begin{equation}\label{multir<p}
\pi_r(\mathfrak{M}_\beta)=\|\beta\|_{r}
\end{equation}  

\item For $p\ge 2$ and $r\ge p$, 
\begin{equation}\label{multir>p}
\pi_r(\mathfrak{M}_\beta)\approx\|\beta\|_{p}
\end{equation} 
\end{enumerate} 
\end{proposition}

{\bf Proof.} 1) Since $\ell^p$ has cotype $2$, an operator from $\ell^p$ to itself is 
$r$-summing operator for some $r\ge1$ if and only if it is $r$-summing for every $r\ge1$. 
Hence it suffices to treat the case $r=2$. We have two different arguments. 
The first one follows from the fact that the composition of two $2$-summing operators is nuclear. 
Here this gives that $\mathfrak{M}_\beta\circ\mathfrak{M}_\beta=\mathfrak{M}_{\beta^2}$ is nuclear. 
It is then easy to conclude that $\beta^2\in\ell^1$.

Another argument uses only the Pietsch domination theorem: we are given $C=\pi_2(\mathfrak{M}_\beta)\ge0$ 
and a probability measure $\nu$ on the unit ball of $(\ell^p)^\ast=\ell^{p'}$ such that, for every $a\in \ell^p$, we have
$$
\dis\|\mathfrak{M}_\beta(a)\|_p^2\le C^2\int_{B_{\ell^{p'}}}\Big|\sum a(n)\alpha(n)\Big|^2\, d\nu(\alpha).
$$

We apply this inequality to 
$$
a_\omega=\dis\sum_{n=0}^Nr_n(\omega)u_n\e_n
$$ 
where $N\ge1$, $(r_n)$ is a Rademacher sequence  and $u$ is a norm $1$ multiplier from $\ell^{p'}$ to $\ell^2$, 
or equivalently, belongs to the unit ball of $\ell^{q}$ with $q=2p/(2-p)$. 
Then we can take the expectation with respect to $\omega$ to obtain via Fubini:  
$$
\dis\|\mathfrak{M}_\beta(u)\|_p^2\le 
C^2\displaystyle\int_{B_{\ell^{p'}}}\|\alpha\cdot u\|_2^2\, d\nu(\alpha)\le 
C^2\displaystyle\int_{B_{\ell^{p'}}}\|\alpha\|^2_{p'}\|u\|^2_q\, d\nu(\alpha)\le C^2.
$$
Hence $\beta$ is actually a multiplier from $\ell^q$ to $\ell^p$ (with norm less than $C$): 
we obtain that $\beta$ belongs to $\ell^2$ (with norm less than $C$).

Conversely, when $\beta\in\ell^2$, we can factorize $\mathfrak{M}_\beta$ through the identity from $\ell^1$ 
to $\ell^2$. 
Indeed, writing $q=2p/(2-p)\in[2,+\infty]$, it suffices to write $\beta=bc$, where $b\in\ell^{p'}$ 
(hence it induces a multiplier from  $\ell^p$ to $\ell^1$), $c\in\ell^{q}$ 
(hence it induces a multiplier from  $\ell^2$ to $\ell^p$) and $\|\beta\|_p=\|b\|_{p'}\cdot\|c\|_{q}$.

2) When $\mathfrak{M}_\beta$ is $r$-summing, its adjoint $\mathfrak{M}_\beta:\ell^{p'}\to\ell^{p'}$ is
order-bounded (see \cite{DJT} p.109), which is equivalent to $\beta\in\ell^{p'}$. 
Conversely, when $\beta\in\ell^{p'}$, we can factorize in an obvious way $\mathfrak{M}_\beta$ 
through the identity from $\ell^1$ to $\ell^2$, which is absolutely summing 
thanks to the Grothendieck's theorem.

3) The sequence $(\e_n)$ is weak-$\ell^r$ (with norm $1$ actually) since $\ell^{r'}\subset\ell^{p}$. 
Hence by definition $\|\beta\|_{r}\le\pi_r(\mathfrak{M}_\beta)$.
Conversely, when $\beta\in \ell^r$, thanks to the fact that $r\le p$, the operator $\mathfrak{M}_\beta$ 
factorizes through the multiplier by $\beta$, viewed from $\ell^\infty$ to $\ell^r$. 
But this multiplier is  $r$-summing with norm less than $\|\beta\|_{r}$.

4) When $\beta\in \ell^p$, the operator $\mathfrak{M}_\beta$ factorizes through 
the multiplier by $\beta$, viewed from $\ell^\infty$ to $\ell^p$, which is  
$p$-summing with norm less than $\|\beta\|_{p}$. 
A fortiori, $\pi_r(\mathfrak{M}_\beta)\le\pi_p(\mathfrak{M}_\beta)\le\|\beta\|_{p}$.

On the other hand, since $\ell^p$ has cotype $p$, 
we deduce from \cite{DJT} p.222, that the operator $\mathfrak{M}_\beta$ is 
$(p,2)$-summing as soon as it is $r$-summing. 
Now, the canonical basis is clearly weak-$\ell^2$ because $p\ge 2$, so $\beta\in\ell^p$.\cqfd
\medskip

Even without  having the full characterization yet, we are now ready to exhibit an example of a composition 
operator on $H^p$, with $p\in(1,2)$, which is order bounded but not absolutely summing (recovering 
the result of Domenig \cite{Do}). 
Indeed, we consider the symbol constructed in 
\cite[proofs of  Th. 4.1 and Lem. 3.7]{LLQR3}
with $\beta\in\dis\big(1,2/p\big)$. With this symbol, the size of the Luecking 
boxes is controlled as follows: for every $n,j$:
$$\dis \lambda_\varphi(R_{n,j})\approx 2^{-n(1+\frac{1}{\beta})}.$$
Applying both (\ref{multi2}) and Th.~\ref{CourSomGen}, 
we get that the operator $J_{\mu_{2^n}}$ has a $2$-summing norm of the 
order of $\dis 2^{(n/2-n/p\beta)}$. Since it is summable, it implies that $J_\mu$ is $2$-summing.

On the other hand, $J_\mu$ is not order bounded since
$$
\int_{\overline{\mathbb D}}\frac{1}{1-|z|}\,d\lambda_\varphi=\int_{\mathbb D}\frac{1}{1-|z|}\,d\lambda_\varphi
\approx\sum_{n,j}2^n\lambda_\varphi(R_{n,j})\approx \sum_n 2^{n(1-\frac{1}{\beta})}=\infty\,. 
$$


\subsection{Some glue-lemmas for summing operators}
\medskip

\begin{lemma}\label{lemreformulationintegrale}
We fix $r\ge1$ and $\eta\in(0,1)$. Let $S:X\to Y$ be an $r$-summing operator. 
There exists a step function $F\colon [0,1]\rightarrow X$ such that 
\begin{itemize}
\item $\dis\int_0^1\big|\chi(F(t))\big|^r\;dt\le1$ for every $\chi$ in the unit ball of $X^\ast$.

\item For every $t\in[0,1]$, $\big\|S(F(t))\big\|\ge\eta\,\pi_r(S)$.
\end{itemize} 
\end{lemma}

{\bf Proof.} By definition of the summing norm $\pi_r(S)$, 
we can find a finite number of vectors $x_1,\ldots,x_n\in X$ such that $S(x_j)\neq0$, 
$$
\sup_{\xi\in B_{X^\ast}}\sum_{1\le j\le n}\big|\xi(x_j)\big|^r\le1 
$$
and 
$$
\Big(\sum_{1\le j\le n}\big\|S(x_j)\big\|^r\Big)^{1/r}\ge\eta\cdot \pi_r(S) . 
$$

Now, choose a mesurable partition $A_1,\ldots,A_n$ of $[0,1]$ such that each $A_j$ has measure 
$$
|A_j|=\frac{\big\|S(x_j)\big\|^r}{\dis\sum_{1\le l\le n}\big\|S(x_l)\big\|^r}\,.
$$
Define 
$$
F=\sum_{j=1}^n\dis\frac{1}{|A_j|^{1/r}}\,x_j\,\ind_{A_j}.
$$
It is now very easy to check that $F$ works.\hfill\cqfd
\medskip


\begin{proposition}\label{sommedirectep>r}
We fix $p\ge r\ge1$. Assume that for every $n\ge1$, there is an $r$-summing operator 
$T_n:X\to Y_n$ such that $\dis\big(\pi_r(T_n)\big)_{n\ge1}\in\ell^r$. 
Then the operator 
\medskip
 
\centerline{\begin{tabular}{rccl}
$T\colon X$ &$\longrightarrow$& $\bigoplus_{\ell^p}Y_n$\cr
\noalign{\smallskip}
$x$& $\longmapsto $& $\dis\big(T_n(x)\big)_{n\ge1} $
\end{tabular}}
\noindent
is an $r$-summing operator and we have
$$
\pi_r(T)\le\Big(\sum_{n\ge1}\pi_r(T_n)^r\Big)^{1/r}.
$$

\end{proposition}

{\bf Proof.} Let us fix vectors  $x_1,\ldots,x_m\in X$  such that
$$
\sup_{\xi\in B_{X^\ast}}\sum_{1\le j\le m}\big|\xi(x_j)\big|^r\le1.
$$
We have
$$
\sum_{1\le j\le m}\big\|T(x_j)\big\|^r= \sum_{1\le j\le m}\Bigg(\sum_{n\ge1}\Big\|T_n (x_j) \Big\|^p\Bigg)^{r/p}.
$$
Since $p\ge r$, we get
$$
\sum_{1\le j\le m}\big\|T(x_j)\big\|^r\le \sum_{1\le j\le m\atop{n\ge1}}\Big\|T_n (x_j) \Big\|^r\le
\sum_{n\ge1}\pi_r(T_n)^r\sup_{\xi\in B_{X^\ast}}\sum_{1\le j\le m}\big|\xi(x_j)\big|^r
$$
\hfill\cqfd


\begin{proposition}\label{sommedirecte2r>2}
We fix $p\ge1$ and $r\ge2$. Assume that for every $n\ge1$, we have a bounded operator $T_n:X_n\to Y_n$. 
We assume that the operator  
\medskip

\centerline{\begin{tabular}{rccl}
$T\colon X=\bigoplus_{\ell^2}X_n$ &$\longrightarrow$& $\bigoplus_{\ell^p}Y_n$\cr
\noalign{\smallskip}
$\dis(x_n)_n$& $\longmapsto $& $\dis\big(T_n(x_n)\big)_{n\ge1} $
\end{tabular}}

\smallskip
\noindent
is an $r$-summing operator.
Then each $T_n$ is $r$-summing and we have
$$
k_r\pi_r(T)\ge\Big(\sum_{n\ge1}\pi_r(T_n)^p\Big)^{1/p}
$$
where $k_r$ is the constant given by the $L^r$-$L^2$-Khinchine inequality.
\end{proposition}

{\bf Proof.} We fix $\eta\in(0,1)$. Thanks to Lemma~\ref{lemreformulationintegrale}, 
we have, for each $n\ge1$, a function $F_n:[0,1]\rightarrow X_n$ such that 
\begin{itemize}
\item $\dis\int_0^1\big|\chi(F_n(t))\big|^r\;dt\le1$,\quad for every $\chi$ in the unit ball of $X_n^\ast$.

\item For every $t\in[0,1]$, $\big\|T_n(F_n(t))\big\|\ge\eta\cdot\pi_r(T_n)$.
\end{itemize} 

Now we consider a Rademacher sequence $(r_n)_{n\ge1}$ (viewed on $[0,1]$) and we can  define the function

\centerline{\begin{tabular}{rccl}
$F\colon$&$[0,1]^2$ &$\longrightarrow$& $\bigoplus_{\ell^2}X_n$\cr
\noalign{\smallskip}
& $\dis(t,\omega)$& $\longmapsto $& $\dis\big(F_n(t)r_n(\omega)\big)_{n\ge1} $
\end{tabular}}

\medskip
\noindent
On one hand, we have for every $\xi$ in the unit ball of 
$X^\ast$: $\xi=(\xi_n)_{n\ge1}$ where $\xi_n$  belongs to $X_n^\ast$ and 
$\dis\sum_{n\ge1}\|\xi_n\|^2\le1$. Thanks to the Khinchine inequality:
$$
\iint_{[0,1]^2}\big|\xi(F(t,\omega))\big|^r\;dtd\omega=
\int_0^1\int_0^1\Big|\dis\sum_{n\ge1}r_n(\omega)\xi_n(F_n(t))
\Big|^r\;d\omega dt\le k_r^r\int_0^1\dis\Big(\sum_{n\ge1}|\xi_n(F_n(t))|^2\Big)^{r/2}\;dt
$$ 
Invoking the triangular inequality in $L^{r/2}$ (recall that $r\ge2$), we get 
$$
\dis\iint_{[0,1]^2}\big|\xi(F(t,\omega))\big|^r\;dtd\omega\le 
k_r^r\dis\bigg(\sum_{n\ge1}\Big(\int_0^1|\xi_n(F_n(t))|^r\;dt\Big)^{2/r}\bigg)^{r/2}\,.
$$

We can write $\xi_n=\|\xi_n\|\chi_n$ where $\chi_n$ belongs to the unit ball of $X_n^\ast$. We obtain
$$
\dis\iint_{[0,1]^2}\big|\xi(F(t,\omega))\big|^r\;dtd\omega\le 
k_r^r\bigg(\sum_{n\ge1}\|\xi_n\|^2\Big(\int_0^1|\chi_n(F_n(t))|^r\;dt\Big)^{2/r}\bigg)^{r/2}
\le k_r^r\bigg(\sum_{n\ge1}\|\xi_n\|^2\bigg)^{r/2}.
$$
Therefore
$$
\bigg(\iint_{[0,1]^2}\big|\xi(F(t,\omega))\big|^r\;dtd\omega\bigg)^{1/r}\le k_r.
$$

On the other hand, 
$$
\iint_{[0,1]^2}\big\|T\big(F(t,\omega)\big)\big\|^r\;dtd\omega=
\iint_{[0,1]^2}\Big(\sum_{n\ge1}\big\|T_n\big(F_n(t)\big)\big\|^p\Big)^{r/p}\;dtd\omega
\ge\eta^r\Big(\sum_{n\ge1}\pi^p_r(T_n)\Big)^{r/p}.
$$
At last, it suffices to observe that 
$$
\bigg(\iint_{[0,1]^2}\big\|T\big(F(t,\omega)\big)\big\|^r\;dtd\omega\bigg)^{1/r}
\le\pi_r(T)\sup_{\xi\in B_{X^\ast}}\bigg(\iint_{[0,1]^2}\big|\xi(F(t,\omega))\big|^r\;dtd\omega\bigg)^{1/r}.
$$
Since $\eta$ is arbitrary, we get the conclusion.\hfill\cqfd

The following result is a variant of the preceding one.

\begin{proposition}\label{sommedirectepr>2}
We fix $q\ge1$ and $r\ge 1$. Assume that for every $n\ge1$, we have a bounded operator $T_n:X_n\to Y_n$.  
We assume that the operator  
\medskip

\centerline{\begin{tabular}{rccl}
$T\colon X=\bigoplus_{\ell^q}X_n$ &$\longrightarrow$& $\bigoplus_{\ell^\infty}Y_n$\cr
\noalign{\smallskip}
$\dis(x_n)_n$& $\longmapsto $& $\dis\big(T_n(x_n)\big)_{n\ge1} $
\end{tabular}}

\smallskip
\noindent
is an $r$-summing operator.
Then each $T_n$ is $r$-summing and we have
$$
\pi_r(T)\ge\Big(\sum_{n\ge1}\pi_r(T_n)^r\Big)^{1/r}\qquad \hbox{when }r\ge q'
$$
and 
$$
\pi_r(T)\ge\Big(\sum_{n\ge1}\pi_r(T_n)^{q'}\Big)^{1/{q'}}\qquad \hbox{when }r\le q'
$$
\end{proposition}

{\bf Proof.} The proof is straightforward:  fix $\eta\in(0,1)$ and for each $n\ge1$, 
choose a finite family of vectors $(x_n^k)$ such that 
$$
\dis\sup_{\|\chi\|_{X^*_n}\le1}\sum_{k}|\chi(x_n^k)|^r\le1\qquad\text{and}\qquad 
\sum_{k}\big\|T_n(x_n^k)\big\|^r\ge\eta\pi_r(T_n)^r.
$$

Now consider a norm one multiplier $a$ from $\ell^{q'}$ to $\ell^r$ and the family of vectors of $X$ 
defined by $v_{n,k}=a_n\big(0,\ldots,x_n^k,0,\ldots\big)$ 
(the {\sl a priori} non zero entry is placed at the $n^{th}$ place). 
For every $\xi$ in the unit ball of $X^\ast$, we have $\xi=(\xi_n)$ with 
$\dis\sum_{n\ge1}\big\|\xi_n\big\|^{q'}\le1$ so
$$
\sum_{n,k}\big|\xi\big(v_{n,k}\big)\big|^r=\sum_{n,k}|a_n|^r\big|\xi_n(x_n^k)\big|^r=
\sum_{n}|a_n|^r\big\|\xi_n\big\|^r\Big(\sum_k\big|\widetilde{\xi_n}\big(x_n^k\big)\big|^r\Big) 
$$
where $\widetilde{\xi_n}$ lies in the unit ball of $X_n^\ast$. Hence
$$
\sum_{n,k}\big|\xi\big(v_{n,k}\big)\big|^r\le \sum_{n}\big|a_n\big|^r\cdot\big\|\xi_n\big\|^r
\le\Big(\sum_{n\ge1}\big\|\xi_n\big\|^{q'}\Big)^{r/q'}\le1\,.
$$
 
On the other hand, by definition of the $r$-summing norm, we have
$$
\pi_r(T)^r\ge\sum_{n,k}\big\|T\big(v_{n,k}\big)\big\|^r=
\sum_{n}\sum_k\big|a_n\big|^r\cdot\big\|T_n\big(x_n^k\big)\big\|^r\ge\sum_n \eta\big|a_n\big|^r\cdot\pi_r(T_n)^r.
$$
Taking the supremum over the norm one multipliers from 
$\ell^{q'}$ to $\ell^r$ and $\eta<1$, we get the result.\hfill\cqfd

For convenience, we state now three corollaries (these are the versions we shall actually use).

\begin{corol}\label{elnuevo}
We fix $p\le2$ and $r\ge 2$. Assume that for every $n\ge1$, we have a bounded operator $T_n:X_n\to Y_n$.  
We assume that the operator  
\medskip

\centerline{\begin{tabular}{rccl}
$T\colon X=\bigoplus_{\ell^p}X_n$ &$\longrightarrow$& $\bigoplus_{\ell^p}Y_n$\cr
\noalign{\smallskip}
$\dis(x_n)_n$& $\longmapsto $& $\dis\big(T_n(x_n)\big)_{n\ge1} $
\end{tabular}}

\smallskip
\noindent
is an $r$-summing operator.
Then each $T_n$ is $r$-summing and we have
$$
k_r\pi_r(T)\ge\Big(\sum_{n\ge1}\pi_r(T_n)^2\Big)^{1/2}\,.
$$
\end{corol}

{\bf Proof.} For every scalar sequence  $(a_n)_n$ which represents a norm one multiplier from 
$\ell^2$ to $\ell^p$, we can apply 
Proposition \ref{sommedirecte2r>2} to the sequence of operators $(a_nT_n)_n$ and we get
$$
k_r\pi_r(T)\ge \Big(\sum_{n\ge1}|a_n|^p\pi_r(T_n)^p\Big)^{1/p}\,.
$$
The result follows taking the supremum over all the norm one multipliers $(a_n)_n$.
\hfill\cqfd

Next result is a direct consequence of Proposition \ref{sommedirectepr>2} for $q=2$, and the fact that
the injection of $\bigoplus_{\ell^p}Y_n$ into $\bigoplus_{\ell^\infty}Y_n$ is  a norm one operator.

\begin{corol}\label{sommedirectepr>2CORr<2}
We fix $p\ge1$ and $r\ge 1$. Assume that for every $n\ge1$, we have a bounded operator $T_n:X_n\to Y_n$.  
We assume that the operator  
\medskip

\centerline{\begin{tabular}{rccl}
$T\colon X=\bigoplus_{\ell^2}X_n$ &$\longrightarrow$& $\bigoplus_{\ell^p}Y_n$\cr
\noalign{\smallskip}
$\dis(x_n)_n$& $\longmapsto $& $\dis\big(T_n(x_n)\big)_{n\ge1} $
\end{tabular}}

\smallskip
\noindent
is an $r$-summing operator.
Then each $T_n$ is $r$-summing and we have
$$
\pi_r(T)\ge\Big(\sum_{n\ge1}\pi_r(T_n)^r\Big)^{1/r}\qquad \hbox{when }r\ge 2
$$
and 
$$
\pi_r(T)\ge\Big(\sum_{n\ge1}\pi_r(T_n)^{2}\Big)^{1/{2}}\qquad \hbox{when }r\le 2
$$
\end{corol}

Using Proposition~\ref{sommedirectep>r} in one direction and Proposition~\ref{sommedirectepr>2} and the norm
one injection of $\bigoplus_{\ell^p}Y_n$ into $\bigoplus_{\ell^\infty}Y_n$ in the opposite direction,
we get the last  corollary:

\begin{corol}\label{sommedirectepr>2COR}
We fix $p\ge2$ and $p\ge r\ge p'$. Assume that for every $n\ge1$, we have a bounded operator $T_n:X_n\to Y_n$.  
Consider the operator 
\medskip

\centerline{\begin{tabular}{rccl}
$T\colon \bigoplus_{\ell^p}X_n$ &$\longrightarrow$& $\bigoplus_{\ell^p}Y_n$\cr
\noalign{\smallskip}
$\dis(x_n)_n$& $\longmapsto $& $\dis\big(T_n(x_n)\big)_{n\ge1} $
\end{tabular}}

\smallskip
\noindent
We have
$$
\pi_r(T)\approx \Big(\sum_{n\ge1}\pi_r(T_n)^r\Big)^{1/r}\,.
$$
\end{corol}


\section{Consequences for Carleson embeddings}

In this section, we will obtain some estimates about summing norms of Carleson embeddings
exploiting the results obtained for the diagonal case (see Section \ref{sectiondiag}) and ``glueing'' the partial 
operators $\widetilde{T}_{p,N}$ or $J_{\mu_N}$. In some cases they provide us with characterizations of
$r$-summingness.
Of course the space $\dis L^p(\D,\mu)$ is the $\ell^p$-sum of the spaces $L^p(\D,\mu_{2^n})$, 
so it is rather easy to glue the range. One of the main difficulty is that, except for $p=2$, $H^p$ 
is not an $\ell^q$-sum 
of a sequence of  spaces $H^p_{N_n}$, whatever may be the value of $q$. 
Nevertheless the Littlewood-Paley 
theorem implies that we can write any $f\in H^p$ as an unconditional sum of $f_j$'s, 
with $f_j\in \widetilde{H^p_{2^j}}$. 

According to the values of $p$ relatively to $2$, we can then exploit the 
type and cotype properties of the spaces $H^p$. In particular:

When $p\le2$, the space $H^p$ has cotype $2$ and type $p$ so we have:
\begin{equation}\label{LPp<2}
\Big(|\hat f(0)|^p+\sum_{j\ge0}\|f_j\|_p^p\Big)^{1/p}\gtrsim\|f\|_p\gtrsim
\Big(|\hat f(0)|^2+\sum_{j\ge0}\|f_j\|_p ^2\Big)^{1/2}.
\end{equation}

When $p\ge2$, the space $H^p$ has cotype $p$ and type $2$ so we have:
\begin{equation}\label{LPp>2}
\Big(|\hat f(0)|^p+\sum_{j\ge0}\|f_j\|_p^p\Big)^{1/p}\lesssim\|f\|_p\lesssim 
\Big(|\hat f(0)|^2+\sum_{j\ge0}\|f_j\|_p ^2\Big)^{1/2}.
\end{equation}

\medskip

Another difficulty we must take care of is the fact that summing up the operators $\dis\widetilde{T}_{p,2^j}$, 
we shall not get precisely the operator $J_\mu$ but its diagonal version. 
More precisely, let us compare the operator $J_\mu$ and the operator 
$$
T_p\colon f\in H^p\longmapsto\sum_{j=0}^{+\infty}f_j\ind_{\Gamma_{2^j}}\in  L^p(\D,\mu).
$$

\begin{lemma}\label{CompJdiag}

Let $p>1$ and $r\ge 1$. 
If the operator $J_\mu$ is $r$-summing, then $T_p$ is $r$-summing and 
$$
\pi_r(T_p)\lesssim\pi_r(J_\mu).
$$
\end{lemma}

{\bf Proof.} We first point out that, for every $f\in H^p$,
$$
\|T_p(f)\|_{ L^p(\D,\mu)}= 
\bigg(\sum_{j=0}^{+\infty}  \Big\|\widetilde{T}_{p,2^j}(f_j)  \Big\|_{ L^p(\D,\mu_{2^j})}  ^p \bigg)^{1/p}\,.
$$

For the proof of the lemma we shall invoke a random perturbation argument like in the proof of 
Proposition~\ref{CourSomGen}:  let $\sigma=(\sigma_0,\sigma_1,\ldots)$ be a sequence of 
Rademacher variables (i.e. independent Bernoulli variables over a probability space $(\Omega,\Prob)$, 
taking their values in $\{\pm1\}$). 
We introduce a function $q_\sigma\in L^\infty(\D,\mu )$ and an operator $R_\sigma\colon H^p\to H^p$ by
$$
q_\sigma=\sum_{j=0}^{+\infty}\sigma_j\ind_{\Gamma_{2^j}}\,\qquad\text{and}\qquad
R_\sigma f=\sum_{j=0}^{+\infty}\sigma_j f_j\quad,\,\hbox{for } f\in H^p\,.
$$
We have $\|q_\sigma\|_\infty=1$, for every $\sigma$, and, thanks to the unconditionality of 
Littlewood-Paley decomposition, 
the operator $R_\sigma\colon H^p\to H^p$ is an isomorphism with bounds not depending on $\sigma$.

It is easy to check that: 
$$
\E_\sigma \Big(\sum_{j=0}^{+\infty}\sigma_jf_j \Big) 
\Big(\sum_{k=0}^{+\infty}\sigma_k\ind_{\Gamma_{2^k}}\Big)=
\sum_{j=0}^{+\infty}f_j\ind_{\Gamma_{2^j}}\,.
$$

So we have for every $f\in   H^p $: 
$\dis\E_\sigma\; q_\sigma J_\mu\circ R_\sigma(f)=T_p(f)$ 
and this implies that for every $r\ge1$, 
$$
\pi_r(T_p)\lesssim\pi_r(J_\mu).
$$
\hfill\cqfd


\medskip

Now we are going to apply the results of the previous 
section to get some estimates of the $r$-summing norms of the Carleson embeddings.

\subsection{Some estimates for $1<p\le2$}

Since $H^p$ and $L^p(\mu)$ have cotype $2$, we recall that, for any $r\ge1$, 
an operator $T:H^p\rightarrow L^p(\mu)$ is $r$-summing if and only if it is $2$-summing.

\begin{proposition}\label{Encadrep<2}
Let  $1<p\le2$ and $\mu$ be a positive measure on the unit disk $\D$. We have
$$
\Bigg[\sum_{n\ge0}\sum_{j=0}^{2^n-1}\Big(2^{n}\mu(R_{n,j})\Big)^{2/p}\Bigg]^{1/2}
\lesssim\pi_2(J_\mu)
\lesssim\Bigg[\sum_{n\ge0}\Big(\sum_{j=0}^{2^n-1}\Big(2^{n}\mu(R_{n,j})\Big)^{2/p}\Big)^{p/2}\Bigg]^{1/p} \,.
$$
\end{proposition}

{\bf Proof.} From (\ref{multi2}) and Theorem~\ref{CourSomGen}, we have the following estimate
for the $2$-summing norm of the operators $J_{\mu_{2^n}}$ and $\tilde T_{p,2^n}$,
$$
\pi_2(J_{\mu_{2^n}})\approx \pi_2(\tilde T_{p,2^n})\approx
\Bigg[\sum_{j=0}^{2^n-1}\Big(2^{n}\mu(R_{n,j})\Big)^{2/p}\Bigg]^{1/2}\,.
$$
We apply (\ref{LPp<2}) and Prop.~\ref{sommedirectep>r} (with $r=p$)
to the sequence of operators $(J_{\mu_{2^n}})_n$ to get the upper estimate for the $p$-summing norm 
(which is equivalent to the $2$-summing norm) for $J_\mu$.

Now let assume that $J_\mu$ is $2$-summing. A fortiori the diagonal operator $T_p$ is $2$-summing with 
$\pi_2(T_p)\lesssim\pi_2(J_\mu)$ (cf Lemma~\ref{CompJdiag}). 
Now applying Corollary~\ref{elnuevo} and \eqref{LPp<2},
we get 
$$
\pi_2\big(T_p\big)\gtrsim\Big(\sum_{n\ge 0}\pi_2\big(\widetilde{T}_{p,2^n}\big)^2\Big)^{1/2}
$$
thanks to the factorization 
$\bigoplus_{\ell^p}\widetilde{H^p_{2^n}} \longrightarrow  H^p \stackrel{T_p}{\longrightarrow} L^p(\D,\mu)
\longrightarrow \bigoplus_{\ell^p}L^p(\D,\mu_{2^n})$ 
which acts on each $\widetilde{H^p_{2^n}}$ as $\widetilde{T}_{p,2^n}$. This gives the minoration.
\hfill\cqfd

\subsection{The case $ r\ge p\ge2$}
\begin{proposition} 
Let  $r\ge p\ge2$ and $\mu$ be a positive measure on the unit disk $\D$. We have
$$
\pi_r(J_\mu)\approx\Big(\int_\D\frac{1}{1-|z|}\;d\mu\Big)^{1/p}.
$$
\end{proposition}

{\bf Proof.} From (\ref{multir>p}) and Theorem~\ref{CourSomGen}, 
we know that $ \pi_r\big(\widetilde{T}_{p,2^n}\big)$ is equivalent to 
the $\ell^p$-sum over $j$ of 
$\big(2^{n}\mu(R_{n,j})\big)^{1/p}$.

Now let assume that $J_\mu$ is $r$-summing. A fortiori the diagonal operator $T_p$ is $r$-summing with 
$\pi_r(T_p)\lesssim\pi_r(J_\mu)$ (cf Lemma~\ref{CompJdiag}). Now 
reasoning as in the second part of the proof of Proposition~\ref{Encadrep<2}, but this time using \eqref{LPp>2}
and applying Proposition~\ref{sommedirecte2r>2}, we get 
$$
\pi_r(J_\mu)\gtrsim
\pi_r\big(T_p\big)\gtrsim\Big(\sum_{n\ge0}\pi_r\big(\widetilde{T}_{p,2^n}\big)^p\Big)^{1/p}\,.
$$
This gives the lower estimate
$$
\pi_r(J_\mu)\gtrsim
\Big(\sum_{n\ge 0}\sum_{j=0}^{2^n-1}2^{n}\mu\big(R_{n,j}\big)\Big)^{1/p}\approx
\Big(\int_\D\frac{1}{1-|z|}\;d\mu\Big)^{1/p}.
$$

The other inequality is clear: the order boundedness easily implies the $p$-summingness, 
hence the $r$-summingness since $r\ge p$.\hfill\cqfd

\subsection{The case $2\le r\le p$}

\begin{proposition}\label{lapobre}
Let  $p\ge r\ge2$ and $\mu$ be a positive measure on the unit disk $\D$. We have
$$
\pi_r(J_\mu)\approx\Bigg[\sum_{n\ge0} \sum_{j=0}^{2^n-1}\Big(2^{n}\mu(R_{n,j})\Big)^{r/p} \Bigg]^{1/r} \,.
$$
\end{proposition}

{\bf Proof.} From (\ref{multir<p}), Theorem~\ref{CourSomGen} and Prop.~\ref{sommedirectep>r}, we get 
in the same way as before
$$
\pi_r\big(J_\mu\big)
\lesssim\Bigg(\sum_{n\ge 0}\sum_{j=0}^{2^n-1}\Big(2^{n}\mu\big(R_{n,j}\big)\Big)^{r/p}\Bigg)^{1/r}.
$$

Now we assume that $J_\mu$ is $r$-summing and we use the same idea than in the previous case, 
replacing Proposition~\ref{sommedirecte2r>2} by Corollary~\ref{sommedirectepr>2CORr<2}.
We get 
$$
\pi_r(J_\mu)\gtrsim \pi_r\big(T_p\big)\ge\Big(\sum_{n\ge0}\pi_r\big(\widetilde{T}_{p,2^n}\big)^r\Big)^{1/r}
\approx\bigg[\sum_{n\ge0} \sum_{j=0}^{2^n-1}\Big(2^{n}\mu(R_{n,j})\Big)^{r/p} \bigg]^{1/r}\,.
$$
This gives the minoration.\cqfd

\subsection{An estimate for $p'\le r\le2$}

We can apply (\ref{multir<p}), Theorem~\ref{CourSomGen} and Proposition~\ref{sommedirectep>r}  to get 

\begin{lemma}\label{lemp'<r}
Let  $p'\le r\le2$ and $\mu$ be a positive measure on the unit disk $\D$. We have
$$
\pi_r\big(J_\mu\big)\le\Big(\sum_{n\ge 0}\pi_r\big(J_{\mu_{2^n}}\big)^r\Big)^{1/r}\approx
\Bigg(\sum_{n\ge0}\sum_{j=0}^{2^n-1}\Big(2^{n}\mu\big(R_{n,j}\big)\Big)^{r/p}\Bigg)^{1/r}.
$$
\end{lemma}

\subsection{Some estimates when $1\le r\le p'\le2$}

In the case $p<2$ we proved in Proposition~\ref{Encadrep<2} that the $r$-summing norm of $J_\mu$ is 
between the $\ell^p$-sum and the $\ell^2$-sum of the $r$-summing norm of the pieces $J_{\mu_{2^n}}$.
Following the same ideas one can prove the following estimates
providing some easy to handle necessary and sufficient conditions. Now the $r$-summing norm of $J_\mu$ turns out to be 
between the $\ell^r$-sum and the $\ell^2$-sum of the $r$-summing 
norm of the pieces $J_{\mu_{2^n}}$.

\begin{proposition}\label{Encadrer<p'}
Let  $p\ge2$, $p'\ge r$ and $\mu$ be a positive measure on the unit disk $\D$. We have
$$
\Bigg[\sum_{n\ge0}\Big(\sum_{j=0}^{2^n-1}\Big(2^{n}\mu(R_{n,j})\Big)^\frac{p'}{p}\Big)^\frac{r}{p'}\Bigg]^{1/r}
\gtrsim\pi_r(J_\mu)\gtrsim
\Bigg[\sum_{n\ge0}\Big(\sum_{j=0}^{2^n-1}\Big(2^{n}\mu(R_{n,j})\Big)^\frac{p'}{p}\Big)^\frac{2}{p'}\Bigg]^{1/2}.
$$
\end{proposition}

Nevertheless this result can also be deduced from the exact charaterization of the $r$-summnig norm
that we will give in Theorem~\ref{THEORr<p'}.


\section{The proof when $p\ge2$}

\medskip

In the previous section we gave some partial results when $p\ge2$. In the case $r\ge 2$, they provided a full characterization for the membership of $J_\mu$ to the class of $r$-summing operators.
In this section we deal with the remaining cases, namely $r\le2$. 
The glueing technics are not sufficient 
in this cases and some new ideas are necessary.

\subsection{The proof in the case $p'< r\le 2$}

We already got an 
inequality in Lemma~\ref{lemp'<r} and 
we shall now prove the reverse inequality. 
\medskip

We are going to define a $H^p$-valued function $\Phi$ which
will play a key role in the proof. We denote $\alpha=1/p'-1/r>0$ and $\rho_n=1-2^{-n}$,
for $n\ge 1$. Define $\Phi\colon [0,+\infty)\to H^p$ by
$$
\Phi (t)(z)=\frac{2^{-n\alpha}}{1-\rho_nz\e^{-2i\pi(t-n)}}\, ,\qquad\hbox{for }z\in \D,\quad t\in[n-1,n), \quad n\ge1\,.
$$
We need the following lemma, where for $f\in H^p$ and $g\in H^{p'}$, the duality bracket
$\langle f,g\rangle$ has the meaning
$$
\langle f,g\rangle = \int_\T f(z) g(\overline z)\,d\lambda(z)\,.
$$
It is known that this makes $H^{p'}$ isomorphic to the dual of $H^p$.

\begin{lemma}\label{LemAp'<r<2}
There exists $C>0$ such that, for every $g\in H^{p'}$:
$$\dis\int_0^{+\infty}\big|\langle \Phi(t),g\rangle\big|^r\, dt\le C\|g\|^r_{ H^{p'}}\; .$$
\end{lemma}

{\bf Proof of Lemma~\ref{LemAp'<r<2}.}
For every $t\in [n-1,n)$, we have by definition of $\Phi$: 
$$
\langle \Phi(t),g\rangle=2^{-n\alpha}g\big(\rho_n\e^{-2i\pi(t-n)}\big)\,.
$$
We get:
$$
\int_0^{+\infty}\big|\langle \Phi(t),g\rangle\big|^r\, dt=
\sum_{n\ge1}2^{-nr\alpha}\int_{\mathbb T}\big|g\big(\rho_n\theta\big)\big|^r\;d\lambda(\theta)=
\int_{\mathbb D}\big|g\big|^r\;d\nu\,,
$$
where the measure $\nu$ is defined for every Borel subset $B$ of $\D$ by 
$$
\nu(B)=\sum_{n\ge1}2^{-nr\alpha}\int_{\mathbb T}\ind_B\big(\rho_n\theta\big)\;d\lambda(\theta).
$$

We point out that the statement of our lemma means exactly that, for the measure $\nu$, we have
$H^{p'}(\D)\subset L^r(\nu)$. In other words (see \cite[Th.~9.4]{D}), we have to prove that 
there exists  constant $C>0$ such that 
\begin{equation}\label{rCarleson}
\nu\bigl({\cal W}(\xi,h)\bigr)\le Ch^{r/p'}\,,\qquad \text{for all $\xi\in\T$ and $0<h<1$,}
\end{equation}
which is usually known as $\nu$ to be a $(r/p')$-Carleson measure.

In order to get that, we fix 
$\xi\in\T$ and $h\in(0,1)$. We can write $2^{-(m+1)}\le h <2^{-m}$ (for some integer $m\ge 0$) and we 
compute
$$
\nu\Big( {\cal W}(\xi,h)\Big)=
\sum_{n\ge1}2^{-nr\alpha}\lambda\Big(\big\{\theta\in\T\,|\;\rho_n\theta\in {\cal W}(\xi,h)\big\}\Big)\le
\sum_{n> m}2^{-nr\alpha} h\lesssim \dis h^{r\alpha+1} = h^{r/p'}.$$
We have \eqref{rCarleson} and the lemma follows.
\hfill\cqfd
\medskip

{\bf Proof of \eqref{p'<r<p}.}
It remains to prove the inequality
\begin{equation}\label{loquefalta}
\Bigg(\sum_{n\ge1}\sum_{j=0}^{2^n-1}\Big(2^{n}\mu\big(R_{n,j}\big)\Big)^{r/p}\Bigg)^{1/r}
\lesssim\pi_r(J_\mu)\,\qquad \text{for $p'<r<2$.} 
\end{equation}

Considering the intervals $I_{n,j}=\big[n-1+2^{-n}j, n-1+2^{-n}(j+1)\big)$, $n\ge 0$, $0\le j < 2^n$, we have, 
for every $t\in I_{n,j}$ 
and every $z\in R_{n,j}$,
$$
\big| \Phi(t)(z)\big| =\Big| \frac{2^{-n\alpha}}{1-\rho_nz\e^{-2i\pi(t-n)}} \Big|\gtrsim
\frac{2^{-n\alpha}} {2^{-n}}=2^{n(1-\alpha)}\,.
$$
Hence, for every $t\in \dis I_{n,j}$:
$$ 
\|\Phi(t)\|_{L^p(\mu)}\gtrsim2^{n(1-\alpha)} \big(\mu (R_{n,j})\big)^{1/p}\,,
$$
and we get
$$ 
\int_0^{+\infty}\|\Phi(t)\|^r_{L^p(\mu)}\,dt=
\sum_{n\ge1}\sum_{j=0}^{2^{n}-1}\int_{I_{n,j}}\|\Phi(t)\|^r_{L^p(\mu)}\,dt
\;\gtrsim \;\sum_{n\ge1}\sum_{j=0}^{2^{n}-1}2^{-n}2^{rn(1-\alpha)}\mu\big( R_{n,j}\big)^{r/p}.
$$
Since $-1+r(1-\alpha)=r/p$, we obtain
$$
\int_0^{+\infty}\|\Phi(t)\|^r_{L^p(\mu)}\,dt\;\gtrsim\;
\sum_{n\ge1}\sum_{j=0}^{2^{n}-1} \Big(2^n\mu\big(R_{n,j}\big)\Big)^{r/p}.
$$

Now, thanks to the described isomorphism between $H^{p'}$ and the dual of $H^{p}$, 
our $r$-summingness assumption on $J_\mu$ implies that 
$$
\pi^r_r(J_\mu)\sup_{\|g\|_{ H^{p'}}\le1}\dis\int_0^{+\infty}\big|\langle \Phi(t),g\rangle\big|^r\, dt
\;\gtrsim \;
\int_0^{+\infty}\|\Phi(t)\|^r_{L^p(\mu)}\,dt.
$$
At last, Lemma~\ref{LemAp'<r<2} gives \eqref{loquefalta}.
\hfill\cqfd

Observe that the given argument to prove \eqref{loquefalta} is still valid for $2\le r\le p$ and could have been 
used in the proof of Proposition~\ref{lapobre}.

\subsection{The proof when $1\le r\le p'\le2$}

Our characterization and the proof in that case are different in nature compared to the preceding ones. 
In Proposition~\ref{Encadrer<p'} we already mentioned some lower and upper estimates 
for the $r$-summing norm.
At the end of this section we will see two examples showing that none of these estimates are 
equivalent to the $r$-summing norm in general.
Hence we cannot obtain a formula looking like in the previous cases.

For proving the characterization \eqref{1<r<p'} stated in section 2  we will use the Poisson integral or Poisson
transform $\cal P$ defined, for any  $f\in L^1(\T,d\lambda)$, by 
$$
{\cal P}[f](z)=\dis\int_\T\frac{1-|z|^2}{|1-\overline{w}z|^2}\, f(w)\;d\lambda(w)\,, \qquad z\in\D.
$$
It is known that ${\cal P}[f]$ is a harmonic function on $\D$. In fact, if $f\in {\mathcal C}(\T)$, 
${\cal P}[f]$ is the solution
to the Dirichlet problem with boundary value $f$; i.e., it represents, for $z\in\D$, the value
of the unique continuous function on $\overline\D$ and harmonic on $\D$ that extends $f$.

We also know that $\cal P$ realizes an isometry from the space 
$$
H^p(\T)=\{f\in L^p(\T) : \hat f(m)=0, \forall m<0\}\,
$$
onto $H^p(\D)$, and so, for any measure $\mu$ on $\D$, the inclusion $J_\mu$ of $H^p(\D)$ in $L^p(\mu)$
is a $r$-summing operator if and only if the Poisson transform ${\cal P}\colon H^p(\T)\to L^p(\mu)$ 
is $r$-summing. 

Moreover, for $1<p<\infty$, $H^p(\T)$ is complemented in  $L^p(\T)$ by the Riesz projection
whose kernel is  
$$
Y=\{f\in L^p(\T) : \hat f(m)=0, \forall m\ge 0\}=\{ Df : f\in H^p(\T)\}\,, 
$$
where $Df(z)=\overline{zf(z)}$, for every $z\in \T$. 
It is plain that 
$$
{\cal P}[Df](z)=\overline{z{\cal P}[f](z)}\,,\qquad\text{for every } z\in\D\,,\quad \text{ and every }f\in H^p(\T)\,.
$$
It follows that ${\cal P}\colon H^p(\T)\to L^p(\mu)$ is $r$-summing if and only if
${\cal P}\colon Y\to L^p(\mu)$ is $r$-summing, and therefore if and only if
${\cal P}\colon L^p(\T)\to L^p(\mu)$ is $r$-summing. As a consequence next theorem implies
(and in fact it is equivalent to) the characterization \eqref{1<r<p'}.

\begin{theorem}\label{THEORr<p'}
Let $p> 2$ and let $\mu$ be a finite measure on the unit disk $\D$. For $\xi\in\T$ define
$$
F_0(\xi)=\big(\mu(\D)\big)^\frac{1}{p}\,,\quad F_n(\xi)=2^{n}\big(\mu({\cal W}(\xi,2^{-n})\big)^\frac{1}{p}\,,\ 
\text{ for $n\ge 1$, \quad and} \quad
F(\xi)=\bigg[\sum_{n\ge1}F_n(\xi)^2\bigg]^{1/2}.
$$
Then the following facts are equivalent:
\begin{enumerate}[(a)]
\item The Poisson transform ${\cal P}$, viewed from 
$L^{p}({\mathbb T})$ to $\dis L^{p}({\mathbb D},\mu)$,  is a $p'$-summing operator.

\item $F\in L^{p'}({\mathbb T})$.

\item  The Poisson transform ${\cal P}$, viewed from $L^{p}({\mathbb T})$ to $L^{p}({\mathbb D},\mu)$,  
is a $1$-summing operator.
\end{enumerate}
Moreover, $\pi_r({\cal P}) \approx\pi_1({\cal P}) \approx\dis \|F\|_{L^{p'}({\mathbb T})}$.

\end{theorem}

In order to prove this theorem, first we have to state and prove several results. 

The following one is probably known from the specialists but it seems not to appear under this form in the literature.
\begin{lemma}\label{LemOB}
Assume $1<p<\infty$ and let $H$ be a Hilbert space and $T\colon L^{p}({\mathbb T})\longrightarrow  H$ 
an operator such that $T^\ast$ is order bounded. Then  $T$ is an absolutely summing operator.
Moreover 
$$
\pi_1(T)\le K_G\Big\|\sup_{x\in B_H}{|T^\ast(x)|(\cdot)}\Big\|_{L^{p'}({\mathbb T})},
$$ 
where $K_G$ is the Grothendieck constant.
\end{lemma}

{\bf Proof.} Let $g\in L^{p'}(\T)$ be a function such that $g\ge |T^*x|$ a.e., for every $x\in B_H$.
Then one can factorize $T^*=B\circ A$, where $A\colon H\to L^\infty$ is defined by
$Ax=T^*x/g$, and $B\colon L^\infty\to L^{p'}$ by $Bf=g\cdot f$. It is clear that $\|A\|\le1$  and $\|B\|=\|g\|_{p'}$.
The result follows since $B^*\colon L^p\to L^1\subset (L^\infty)^*$ has norm $\|g\|_{p'}$, $T=A^*|_{L^1}\circ B^*$, and $A^*\colon L^1\to H$ 
is $1$-summing thanks to Grothendieck theorem.

\ \hfill\cqfd
\medskip

The following lemma is a substitute to the  glue lemmas used for the other cases ($r>p'$). 
Instead of gluing some absolutely summing (partial) operators, 
we are going to glue some order bounded (partial) operators.

\begin{lemma}\label{LemOBglue}
Let  $(U,\nu)$ and $(\Omega,\mu)$ be two 
measure spaces, $p>2$, and 
$T\colon L^{p}(U,\nu)\to  L^{p}(\Omega,\mu)$ a bounded operator.

We assume that there exist  sequences of pairwise disjoint measurable subsets $\Omega_m$ of $\Omega$,  
of Hilbert spaces $H_m$,
of operators $A_m\colon L^{p}(U,\nu)\longrightarrow  H_m$, 
of contractions $B_m\colon H_m\longrightarrow  L^{p}(\Omega,\mu)$ and 
of functions $F_m\in L^{p'}(U,\nu)$ with the following properties:

\begin{itemize}
\item $\Omega=\bigcup\Omega_m$.

\item For every $f\in L^{p}(U,\nu)$, we have $\dis B_m\circ A_m (f)=\ind_{\Omega_m}T(f)$.

\item For every $x$ in the unit ball of $H_m$, we have $\dis  |A_m^\ast (x)|\le F_m\;$ $\nu$-a.e. on $U$.

\item The function $F=\dis\big(\sum F_m^2\big)^\frac{1}{2}$ belongs to $\dis L^{p'}(U,\nu)$. 
\end{itemize}
Then,  the operator $T$ is absolutely summning and $\pi_1(T)\le K_G\|F\|_{ L^{p'}(U,\nu)}$.
\end{lemma}

{\bf Proof.} It is natural to consider the Hilbert space $H=\oplus_{\ell^2}H_m$ and the (diagonal) operator 
$B\colon H\to L^p(\Omega)$ defined by $B\big((x_m)\big)=\sum_m \ind_{\Omega_m}\cdot B_m(x_m)$. 
This operator is clearly bounded and actually it is a contraction since
$$
\big\|B\big((x_m)\big)\big\|_p=\Big(\sum_m\big\|\ind_{\Omega_m}\cdot B_m(x_m)\big\|_p^p\Big)^\frac{1}{p}
\le\Big(\sum_m\big\|B_m(x_m)\big\|_p^2\Big)^\frac{1}{2}
\le\Big(\sum_m\big\|x_m\big\|_{H_m}^2\Big)^\frac{1}{2}
=\|x\|_{H}.
$$

On the other hand, the operator $A: L^p(U,\nu)\longrightarrow H $ defined by $A(f)=\big(A_m(f)\big)_m$ 
is bounded as well since its adjoint is so. Indeed, we claim that $A^\ast$ is even order bounded: 
for every $y=(y_m)\in H$, we have $A^\ast(y)=\sum_m A^\ast_m(y_m)$ and for a.e. $\xi\in U$, we have
$$
|A^\ast(y)(\xi)|\le\dis\sum_m| A^\ast_m(y_m)(\xi)|\le \sum_m\|y_m\|F_m(\xi)\le
\Big(\sum_m\big\|y_m\big\|_{H_m}^2\Big)^\frac{1}{2}
\Big(\sum_m F_m(\xi)^2\Big)^\frac{1}{2}
=\|y\| F(\xi)\,,
$$
which was our claim, since $F\in L^{p'}(U,\nu)$ by hypothesis. 

It is easy to check that $T=B\circ A$, and, by Lemma~\ref{LemOB}, we get that $A$ is absolutely summing and 
$$
\pi_1(T)\le\|B\|\pi_1(A)\le K_G\|F\|_{ L^{p'}(U,\nu)}.
$$
\hfill\cqfd

We need now some more specific estimates on the Poisson transform. Recall that, for $a\in \D$, we have
$P_a(\xi)=(1-|a|^2)/|1-\overline a \xi|^2$, for $\xi\in\T$.

\begin{proposition}\label{propPoisson}
Let $\mu$ be a finite measure on the unit disk $\D$, $p>2$ and $a\in\D$. 
Assume $E$ is a Borel subset of $D\big(a,\frac{1-|a|}{2}\big)$. Then
there exist a Hilbert space $H$, an operator  $A\colon L^{p}({\mathbb T}) \to H$ and a contraction 
$B\colon H\to  L^{p}({\mathbb D},\mu)$ such that

\begin{itemize} 

\item For every $f\in L^{p}({\mathbb T})$, we have $ B\circ A(f)=\ind_{E}{\cal P}(f)$.

\item For every $f$ in the unit ball of $H$, we have $|A^\ast (f)(\xi)|\le100 \big(\mu(E)\big)^\frac{1}{p}P_a(\xi)$,
for almost every $\xi\in\T$.
  
\end{itemize}
\end{proposition}

{\bf Proof.} Let $\gamma$ be the boundary of the disk $D\big(a,\frac{3(1-|a|)}{4}\big)$. Our Hilbert space 
$H$ will be $L^2\big(\gamma,m_\gamma\big)$ where $m_\gamma$ is the normalized arc length measure on
$\gamma$. Let $B_1\colon H\longrightarrow L^p(E,\mu)$ be the Poisson operator associated to the domain   
$D\big(a,\frac{3(1-|a|)}{4}\big)$.

For the classical Poisson transform $\cal P$, we have, if $|z|\le 2/3$,
\begin{equation}\label{PoissonEval2}
 |{\cal P}(f)(z)|\le\sqrt3  \|f\|_{L^2(\mathbb T)}\,,\qquad \text{for every $f\in L^2(\mathbb T)$.}
\end{equation}
Indeed, if $|z|\le 2/3$, it is easy to check that $\|P_z\|^2_{L^2(\mathbb T)}\le 13/5\le3$. 
The translation of \eqref{PoissonEval2} to our setting yields
$$
|B_1(g)(z)|\le\sqrt 3 \|g\|_H\,,\qquad\text{for every $g\in H$ and every $z\in E$.}
$$
Hence
$$
\|B_1(g)\|_{L^p(E,\mu)}\le\sqrt3\|g\|_H\big(\mu(E)\big)^{1/p}\,,\qquad\text{for every $g\in H$.}
$$
Therefore defining $Bg=(\ind_{E}\cdot B_1g)/\sqrt3\big(\mu(E)\big)^{1/p}$,  $g\in H$, we have a contraction
from  $H$ to $L^{p}({\mathbb D},\mu)$.

\smallskip

Now we define the operator $A:L^{p}({\mathbb T}) \to H $ as the (classical) 
Poisson integral (restricted to $\gamma$) up to a constant. More precisely 
$A(f)(w)=\sqrt3\big(\mu(E)\big)^{1/p}{\cal P}[f](w)$, for every  $w\in\gamma$
and $f\in L^p(\T)$. Clearly we have 
$ B\circ A(f)=\ind_{E}{\cal P}(f)$, for all  $f\in L^{p}({\mathbb T})$.

It is easy to check that the adjoint of $A$ is given by
$$
A^\ast(g)(\xi)=\sqrt3\big(\mu(E)\big)^{1/p}\int_\gamma P_w(\xi)\cdot g(w)\;dm_\gamma(w) \,,\qquad
g\in H.
$$
When $w\in\gamma$ and $\xi\in\T$, we have 
$$
1-|w|\le \frac{7}{4}( 1-|a|)\,,\qquad 
|w-\xi|\ge|a-\xi|-\frac{3}{4}(1-|a|)\ge\frac{1}{4}|a-\xi|
$$
and
$$
P_w(\xi)=\frac{ 1-|w|^2}{|w-\xi|^2}\le\frac{32( 1-|w|)}{|a-\xi|^2}\le 56 P_a(\xi).
$$
We then obtain the order boundedness of $A^\ast$. 
Indeed, for $\xi\in\T$ and $g\in H$, we have
$$
|A^\ast(g)(\xi)|=\dis\sqrt3 \big(\mu(E)\big)^{1/p}\int_\gamma P_w(\xi)|g(w)|\;dm_\gamma
\le 56\sqrt3\big(\mu(E)\big)^{1/p}P_a(\xi)\|g\|_{L^1(\gamma,m_\gamma)}.
$$
Therefore
$$
\sup_{g\in B_H}|A^\ast(g)(\xi)|\le100\big(\mu(E)\big)^{1/p}P_a(\xi).
$$\hfill\cqfd

For the next result, 
we need some notations: 
for $n\ge1$, we define ${\cal E}_n=\big\{a_{n,1},\ldots,a_{n,m_n}\big\}$ as a maximal 
$2^{-(n+2)}$-net in the dyadic corona $\Gamma_n=\big\{z\in\D\, ;\; 1-2^{-n}\le |z|<1-2^{-(n+1)}\big\}$. 
We also define ${\cal E}_0=\{0\}$ and $E_0=\frac{1}{2}\D$.
We will assume that, for $n\ge 1$,
$\{E_{n,j}\}_{1\le j\le m_n}$ is a family of pairwise disjoint subsets of 
$D\big(a_{n,j},2^{-(n+2)}\big)\cap\Gamma_n$, whose union is the corona $\Gamma_n$.

\begin{proposition}\label{propG}
Let $p>2$ and let $\mu$ be a finite measure on the unit disk $\D$.  For $\xi\in\T$ define
$$
G_0(\xi)=\big(\mu\big(\frac{1}{2}\D\big)\big)^\frac{1}{p} 
\,,\qquad G_n(\xi)=\Big(\sum_{1\le j\le m_n}\big(\mu(E_{n,j})\big)^\frac{2}{p}\big(P_{a_{n,j}}\big)^2\Big)^\frac{1}{2}
\,, \text{ for $n\ge 1$, }
$$
and 
$$
G(\xi)=\big(\sum_{n\ge1}G_n(\xi)^2\big)^{1/2}.
$$
If $G\in  L^{p'}({\mathbb T})$, then ${\cal P}\colon L^{p}({\mathbb T})\to L^{p}({\mathbb D},\mu)$ is
absolutely summing and moreover
$$
\pi_1({\cal P})\lesssim \|G\|_{p'}\, .
$$
\end{proposition}

{\bf Proof.} Observe that 
$D\big(a_{n,j},2^{-(n+2)}\big)\cap\Gamma_n\subset D\big(a_{n,j},\frac{1-|a_{n,j}|}{2}\big)$. Then,
for each fixed $(n,j)$, Prop.~\ref{propPoisson} applies with $E=E_{n,j}$ and $a=a_{n,j}$. 
Lemma~\ref{LemOBglue} can be applied to the (countable) collections of the sets $\Omega_m=E_{n,j}$ and 
functions $F_m=100 \big(\mu(E_{n,j})\big)^\frac{1}{p}P_{a_{n,j}}$ to get the result.
\hfill\cqfd

At last, we can prove Theorem~\ref{THEORr<p'}. 
One of the key point will be to show that the function $G$ of Proposition~\ref{propG} 
is equivalent to the function $F$ of Theorem~\ref{THEORr<p'}.
\medskip

{\bf Proof of Theorem~\ref{THEORr<p'}.}   

{\sl (c) $\Rightarrow$ (a)} is obvious.
\medskip

{\sl (a) $\Rightarrow$ (b)} 
We are going to use some dyadic test functions: 
let us consider for $\xi\in\T$ and $z\in\overline{\D}$, the sequence of functions
$$
K_n(\xi,z)=\dis\sum_{j=2^n}^{2^{n+1}-1}(z\overline{\xi})^j=
(z\overline{\xi})^{2^n}\frac{1-(z\overline{\xi})^{2^n}}{1-z\overline{\xi}}.
$$

Let $r>2$ be such that $\frac{1}{p'}=\frac{1}{r}+\frac{1}{2}$. 
With such a choice of $r$, since $\ell^2$ is isometric to the space of multipliers 
from $\ell^r$ to $\ell^{p'}$, we can choose, for each $\xi\in\T$, a sequence of positive functions 
$g_n(\xi)$ (measurable as functions of $\xi$) such that 
$$
\sum_{n\ge0}(g_n(\xi))^r=1\,,\qquad\text{and}\qquad \dis\sum_{n\ge0}(g_n(\xi))^{p'}|F_n(\xi)|^{p'}=(F(\xi))^{p'}.
$$

Now, we test the $p'$-summingness of $\cal P$ on the random ($ L^{p}({\mathbb T})$ valued) 
function 
$$
(n,\xi)\in\N\times\T\mapsto K_n(\xi,\cdot)g_n(\xi)\,.
$$
We get
$$
\sum_{n\ge0}\int_{\mathbb T}\Big\|K_n(\xi,\cdot)g_n(\xi)\Big\|_{L^{p}({\mathbb D},\mu)}^{p'}\;d\lambda(\xi)
\le\pi_{p'}^{p'}({\cal P})\sup_{h\in B_{L^{p'}({\mathbb T})}}
\sum_{n\ge0}\int_{\mathbb T}
\Big|\int_{\mathbb T}h(z)g_n(\xi)\overline{K_n(\xi,z)}\;d\lambda(z)\Big|^{p'}\;d\lambda(\xi)\,.
$$
But
$\int_{\mathbb T}h(z)g_n(\xi)\overline{K_n(\xi,z)}\;d\lambda(z)=g_n(\xi)Q_n(h)(\xi)$ 
where $Q_n(h)$ is the $n$'th dyadic projection of $h$, i.e. 
$$
Q_n(h)(\xi)=\sum_{j=2^n}^{2^{n+1}-1}\hat h(j)\xi^j.
$$
Now we use the H\"older inequality to majorize $\sum_{n\ge0}\big|g_n(\xi)Q_n(h)(\xi)\big|^{p'}$ by 
$\big(\sum_{n\ge0}\big|Q_n(h)(\xi)\big|^2\big)^{p'/2}$.
This yields
$$
\sum_{n\ge0}\int_{\mathbb T}\Big\|K_n(\xi,.)g_n(\xi)\Big\|_{L^{p}({\mathbb D},\mu)}^{p'}\;d\lambda(\xi)
\le\pi_{p'}^{p'}({\cal P}) \sup_{\|h\|_{L^{p'}({\mathbb T})}\le1}
\int_{\mathbb T}\Big(\sum_{n\ge0}\big|Q_n(h)(\xi)\big|^2\Big)^{p'/2} \;d\lambda(\xi),
$$
and, thanks to the Littlewood-Paley inequality (introducing a constant $C$ depending only on $p'$),
$$
\sum_{n\ge0}\int_{\mathbb T}\Big\|K_n(\xi,.)g_n(\xi)\Big\|_{L^{p}({\mathbb D},\mu)}^{p'}\;d\lambda(\xi)
\le C\pi_{p'}^{p'}({\cal P})\sup_{\|h\|_{L^{p'}({\mathbb T})}\le1}\|h\|_{L^{p'}({\mathbb T})}^{p'}
\le C\pi_{p'}^{p'}({\cal P}).
$$

On the other hand, we notice that $|K_n(\xi,z)|\gtrsim 2^{n}$ when $z\in {\cal W}(\xi,2^{-n})$. 
So the left hand term is bounded below by
$$
\int_{\mathbb T}\sum_{n\ge0}\Big\|K_n(\xi,.)g_n(\xi)\Big\|_{L^{p}({\mathbb D},\mu)}^{p'}\;d\lambda(\xi)\gtrsim 
\int_{\mathbb T}\sum_{n\ge0}( g_n(\xi))^{p'}
\Big(2^n\mu \big({\cal W}(\xi,2^{-n})\big)^\frac{1}{p}\Big) ^{p'}\;d\lambda(\xi) .$$
Then, by the choice of the $g_n$'s, we have 
$$
\dis\sum_{n\ge0}| g_n(\xi)|^{p'}\Big(2^n\mu \big({\cal W}(\xi,2^{-n})\big)^\frac{1}{p}\Big) ^{p'}=
\sum_{n\ge0}( g_n(\xi))^{p'}(F_n(\xi))^{p'}=F^{p'}(\xi).
$$
Therefore we have $F\in L^{p'}({\mathbb T})$ and $(b)$ is proved.

\medskip

{\sl (b) $\Rightarrow$ (c)}  Let us assume that $F\in L^{p'}({\mathbb T})$. 
Assertion $(c)$ will be proved by Prop.~\ref{propG} as soon as we show that $G\lesssim F$ a.e. on $\T$.
Obviously $G_0\le F_0$. 

We fix $n\ge1$ and $\xi\in\T$.
For $l\in\{1,\ldots,n\}$, we set $I_l=\big\{ j\in\{1,\dots, m_n\}\,|\; E_{n,j}\subset W\big(\xi,2^{-l}\big)\big\}$. 
Clearly, 
$$
I_n\subset I_{n-1}\subset\ldots\subset I_0=\{1,\ldots,m_n\}.
$$

It is easy to see that $m_n\approx 2^n$ and that, for every $\xi\in\T$ and $l\in\{1,\ldots,n\}$, 
the cardinality of ${\cal E}_n\cap {\cal W}(\xi,2^{-l})$ is less than $2^{n-l}$ (up to a universal constant).
Now we make these sets disjoints: let $J_n=I_n$ and $J_l=I_l\setminus I_{l+1}$ for $0\le l<n$. 
We have in particular $|J_l|\lesssim2^{n-l}$. 

Moreover, for any $l\in\{1,\ldots,n\} $ we have the following estimate using the H\"older inequality:
$$
\sum_{j\in J_l}\big(\mu(E_{n,j})\big)^\frac{2}{p}\le 
|J_l|^{1-\frac{2}{p}}\Big(\sum_{j\in J_l}\big(\mu(E_{n,j})\Big)^\frac{2}{p}$$
and since $\sum_{j\in J_l}\big(\mu(E_{n,j})=\mu\big(\bigcup_{j\in J_l} E_{n,j}\big)\le
\mu\big({\cal W}(\xi,2^{-l}\big)$, we get by the definition of $F_l$:
$$
\sum_{j\in J_l}\big(\mu(E_{n,j})\big)^\frac{2}{p}\lesssim 2^{(n-l)(1-\frac{2}{p})}2^{-2l}F_l^2(\xi).
$$
When $l=0$, it is actually still valid.

Now, by definition, for any $j\in J_l$, with $l<n$, we have $ E_{n,j}\not\subset  W\big(\xi,2^{-(l+1)}\big)$, 
so there exists some $z\in E_{n,j}$ such that either $1-|z|>2^{-(l+1)}$ or $|\arg(z\overline{\xi})|\ge2^{-(l+1)}$. 
But, since $z\in E_{n,j}\subset\Gamma_n$, we have $1-|z|\le2^{-n}\le2^{-(l+1)}$ because $n\ge l+1$. 
Therefore, we have  $|\arg(z\overline{\xi})|\ge2^{-(l+1)}$. We obtain:
$$
|z-\xi|\ge|z-|z|\xi|\ge|z|.\frac{2}{\pi} |\arg(z\overline{\xi})|\ge\frac{2^{-l}}{2\pi},
$$
since $|z|\ge\frac{1}{2}$.
But $z\in E_{n,j}\subset D\big(a_{n,j},2^{-(n+2)}\big)$, and $2^{-(n+2)}\le2^{-l}/8$, because $l<n$. 
Therefore $ |\xi-a_{n,j}|\ge|z-\xi|-|z-a_{n,j}|\gtrsim2^{-l}$.
The Poisson kernel is then majorized: 
$$
P_{a_{n,j}}(\xi)=\frac{1-|a_{n,j}|^2}{ |\xi-a_{n,j}|^2}\lesssim \frac{2^{-n}}{2^{-2l}}
$$ 
for every $j\in J_l$, when $0\le l<n$. Actually, the same 
estimate is valid when $l=n$.

\medskip

We are now in position to conclude: 
$$ 
G_n^2(\xi)=\sum_{j=0}^{m_n} \big(\mu(E_{n,j})\big)^\frac{2}{p}P^2_{a_{n,j}}(\xi)=
\sum_{l=0}^{n} \sum_{j\in J_l}\big(\mu(E_{n,j})\big)^\frac{2}{p}P^2_{a_{n,j}}(\xi).
$$
Using the preceding estimates, we get
$$
G_n^2(\xi)\lesssim \sum_{l=0}^{n}\frac{2^{-2n}}{2^{-4l}}.2^{(n-l)(1-\frac{2}{p})}2^{-2l}F_l^2(\xi)
=\sum_{l=0}^{n}a^{(n-l)}F_l^2(\xi)$$
where $a=\dis2^{-(\frac{2}{p}+1)}<1$. 
At last, we have 
$$
G^2(\xi)\lesssim\sum_{n\ge0}\sum_{l=0}^{n}a^{(n-l)}F_l^2(\xi)=
\sum_{l\ge0}\sum_{l=n}^{+\infty}a^{(n-l)}F_l^2(\xi)\lesssim\sum_{l\ge0}F_l^2(\xi)=F^2(\xi).
$$\hfill\cqfd 
\bigskip

It is natural to wonder whether there is a ``continuity" in $r=p'$ in our characterizations. More precisely, do we have 
$$
\pi_{p'}(J_\mu)\approx
\Bigg[\sum_{n\ge 0}\;\sum_{0\le j<2^n}\Big(2^{n}\mu(R_{n,j})\Big)^{{p'}/p}\Bigg]^{1/{p'}}\; ?
$$
The answer is negative as soon as $p>2$ and we shall even prove that it is negative for some 
pull back measure associated to a symbol $\varphi$. In other words, it is false even in the class 
of composition operators. 

\begin{ex}\label{Noupper}
For every $p>2$
there exists a symbol $\varphi\colon\D\to\D$ such that $C_\varphi \colon H^p\to H^p$ is $1$-summing but
$$
\sum_{n\ge 0}\;\sum_{0\le j<2^n}\Big(2^{n}\lambda_\varphi(R_{n,j})\Big)^{{p'}/p}=+\infty\,.
$$
\end{ex}

{\bf Proof.} We shall choose a rotation invariant probability measure on $\D$, i.e. 
$\mu(\theta B)=\mu(B)$, for every Borel $B$ and every
$\theta\in\T$, satisfying the condition
\begin{equation}\label{Bishop}
\int_{\mathbb D}\log(1/|z|)\,d\mu(z)<+\infty.
\end{equation}
Then, thanks to \cite[Th. 1.1]{B} there exists a symbol $\varphi:\D\to\D$ with $\varphi(0)=0$ whose associated
pullback measure is $\lambda_\varphi=\mu$. To obtain condition \eqref{Bishop} it is enough that
$\mu$ to be null in a neighbourhood of $0$.

So our purpose is reduced to  find a  probability measure $\sigma$ on $(0,1)$, the measure satisfying  
$\sigma\big((a,b)\big)=\mu(\{z\in\D|\; a<|z|<b\})$, which allows to describe $\mu$ in the following way:
$$ 
\mu(B)=\int_0^1\int_{\mathbb T}\ind_B(rz)\;d\lambda(z)\;d\sigma(r)\,,\qquad\text{for every Borel set $B$.}
$$
In particular, for $\xi\in \T$, $n\ge 0$, and $0\le j <2^n$, we have
$$
\mu\big({\cal W}(\xi, 2^{-n})\big)=\frac{2^{-n}}{\pi}\sigma\big([1-2^{-n},1)\big)\,,\qquad\text{and}\qquad
\mu(R_{n,j})=2^{-n}\sigma\big([1-2^{-n},1-2^{-n-1})\big)\,.
$$
Then, if we call $x_n=\sigma\big([1-2^{-n},1-2^{-n-1})\big)$, we need
\begin{equation}\label{infinito}
+\infty=\sum_{n\ge 0}\;\sum_{0\le j<2^n}\big(2^{n}\mu(R_{n,j})\big)^{{p'}/p}=
\sum_{n\ge 0} 2^{n}x_n^{p'/p}\,.
\end{equation}

On the other hand, if we call $y_n= \sigma\big([1-2^{-n},1)\big)$, we have, for every $\xi\in\T$,
$$
\pi^{1/p} F(\xi)=\Big[\sum_{n\ge0} 2^{2n}(2^{-n}y_n)^{2/p}\Big]^{1/2}
=\Big[\sum_{n\ge0} 2^{2n/p'}y_n^{2/p}\Big]^{1/2}
$$
Consequently, as $F$ is a constant function, if we have 
\begin{equation}\label{finito}
\sum_{n\ge0} 2^{2n/p'}y_n^{2/p} < +\infty,
\end{equation}
by Theorem~\ref{THEORr<p'}, the operator $J_\mu$ (and then $C_\varphi$) is $1$-summing.

We are going to choose the probability $\sigma $, given a positive decreasing sequence $(\alpha_n)$
and putting $z_n=1-2^{-n}$, by the sum
$$
\sigma=\sum_{n\ge 1} 2^{-np/p'} \alpha_n \delta_{z_n}\,.
$$
Then $x_n=2^{-np/p'}\alpha_n$, and it is clear that 
$y_n=\sum_{m\ge n} x_m\le C2^{-np/p'}\alpha_n$.
Therefore
$$
\sum_{n\ge 1} 2^{n}x_n^{p'/p}=\sum_n \alpha_n^{p'/p}\,,\qquad\text{and}\qquad
\sum_{n\ge1} 2^{2n/p'}y_n^{2/p} \le C^{2/p} \sum_{n\ge1}\alpha_n^{2/p}\,.
$$
So, in order to have \eqref{infinito} and  \eqref{finito} it suffices to choose the sequence
$(\alpha_n)_n$ in $\ell^{2/p}\setminus \ell^{p'/p}$, and this is possible since $2>p'$.
\hfill\cqfd

The previous example also shows that the upper estimate in Proposition~\ref{Encadrer<p'} is not equivalent
to the $r$-summing norm in general. In the next example we show that the same happens with
the lower estimate given in Proposition~\ref{Encadrer<p'}.

\begin{ex}\label{Nolower}
For every $p>2$
there exists a finite measure $\mu$ on $\D$ such that
\begin{equation}\label{aquiarriba}
\Bigg[\sum_{n\ge0}\Big(\sum_{j=0}^{2^n-1}
\Big(2^{n}\mu(R_{n,j})\Big)^\frac{p'}{p}\Big)^\frac{2}{p'}\Bigg]^{1/2}<+\infty ,
\end{equation}
but $J_\mu\colon H^p\to L^p(\mu)$ is not $p'$-summing.
\end{ex}

{\bf Proof.} The measure $\mu$ is going to be supported in a radius of $\D$, concretely in the segment $[0,1]$.
Moreover $\mu$ will be of the form
$$
\mu=\sum_{n\ge 1} \alpha_n \delta_{z_n}\,,
$$
where $\alpha_n>0$, $\sum_n\alpha_n <+\infty$, and $z_n=1-2^{-n}$. In this case we will have \eqref{aquiarriba}
as soon as 
\begin{equation}\label{estopasa}
+\infty >  \sum_{n\ge0}
\big(2^{n}\mu(R_{n,0})\big)^{2/p} = \sum_{n=1}^\infty\big(2^{n}\alpha_n\big)^{2/p} .
\end{equation}

To find $\alpha_n$'s  for $J_\mu$ not to be $p'$-summing, we could use the characterization
in Theorem~\ref{THEORr<p'}, but we are going to provide a different argument.
Consider, for the points $z_n$, their reproducing kernels 
$$
K_{z_n}(w)=\frac{1}{1-w\overline{z_n}}\,, \qquad w\in\D,
$$
and define $f_n=2^{-n/p'}K_{z_n}$. We claim that the sequence $(f_n)$ is weakly $p'$-summable in $H^p$.
Then, if $J_\mu$ is $p'$-summing, we have
\begin{equation}\label{estonopasa}
+\infty >  \sum_{n\ge1} \|f_n\|_{L^p(\mu)}^{p'} \ge \sum_{n\ge 1} |f_n(z_n)|^{p'}\alpha_n^{p'/p}
\gtrsim  \sum_{n\ge 1} (2^{-n/p'} 2^n)^{p'}\alpha_n^{p'/p} = \sum_{n\ge 1} (2^n\alpha_n)^{p'/p} .
\end{equation}
Then it is clear that we can choose $\alpha_n$'s satisfying \eqref{estopasa} and not  \eqref{estonopasa}
since $p'<2$.

For proving our claim, take into account that $\langle g,K_{z_n} \rangle= g(z_n)$, for every $g\in H^{p'}$.
Then
$$
\sum_{n\ge 1} |\langle g,f_n \rangle| ^{p'}=\sum_{n\ge 1} 2^{-n}|g(z_n)|^{p'} =\int_{\mathbb D} |g|^{p'}\, d\nu\;,
$$ 
where $\nu$ is the measure $\nu=\sum_{n\ge1} 2^{-n}\delta_{z_n}$. It is easy to see that $\nu$ is a Carleson
measure and then there exists $C$, such that 
$$
\int_{\mathbb D} |g|^{p'} \, d\nu \le C^{p'} \|g\|^{p'}_{H^{p'}}\,,\qquad \text{for every $g\in H^{p'}$.}
$$
The claim and the example follow.
\hfill\cqfd

\section{The case $1<p\le2$}

Before giving the main results of this section let us state a proposition which yields in particular the
equivalence $\|\Phi\|_{2/p}^{1/p}\approx \|\Psi\|_{2/p}^{1/p}$ in \eqref{p<2}.

\begin{proposition}\label{equivInt}
Let $\nu$ be a positive finite measure on the unit disk $\D$, $\gamma\ge1$ and $\eta>1$.
We have the equivalences
$$
\int_\T \biggl( \int_\D \frac{1}{ |1-\overline{w}z|^\eta}\;d\nu(z)  \, \biggr)^{\gamma} \, d\lambda(w)\Bigg)^{1/\gamma}\approx
\Bigg(\int_\T \biggl( \int_{\Sigma_w} \frac{1}{|1-\overline{w}z|^\eta}\;d\nu(z)  \, \biggr)^{\gamma} 
\, d\lambda(w)\Bigg)^{1/\gamma}
$$
$$
\approx \Bigg(\int_\T \biggl( \int_{\Sigma_w} \frac{1}{(1-|z|^2)^\eta}\;d\nu(z)  \, \biggr)^{\gamma} 
\, d\lambda(w)\Bigg)^{1/\gamma},
$$
with constants depending only on $\gamma$ and $\eta$.
\end{proposition}

Let us remark that when $\gamma=1$, we have by Fubini and \cite[Th. 1.7]{HKZ}:
\begin{equation}\label{hkz}
\int_\T \int_\D \frac{1}{ |1-\overline{w}z|^\eta}\;d\nu(z) \, d\lambda(w)\approx
\int_\D \frac{d\nu}{ (1-|z|^2)^{\eta-1} }\;.
\end{equation}

{\bf Proof of Prop~\ref{equivInt}}.

Let $w\in\T$. Of course, we always have $\dis1-|z|^2\le2(1-|z|)\le2|1-\overline{w}z|$ but we point out that these 
quantities are actually equivalent on the Stolz domain  $\Sigma_w$: for every $z\in\Sigma_w$, 
we have $\dis1-|z|^2\approx|w-z|=|1-\overline{w}z|$ (up to numerical constants). 
This proves that the two last quantities in the statement are equivalent. 
We also get obviously that the first integral is greater (up to constants) than the third one.

Now let us prove the converse. We wish to prove that $ B_\nu\lesssim A_\nu$ where 
$$
A_\nu=
\bigg(\int_\T \biggl( \int_{\Sigma_w} \frac{1}{(1-|z|^2)^\eta}\;d\nu(z)  \, \biggr)^{\gamma} \, 
d\lambda(w)\bigg)^{1/\gamma}\hbox{and}\quad 
B_\nu= \bigg(\int_\T \biggl( \int_\D \frac{1}{ |1-\overline{w}z|^\eta}\;d\nu(z)  \, \biggr)^{\gamma}\, d\lambda(w)\bigg)^{1/\gamma}
$$
First we linearize ($\gamma'$ is the conjugate exponent of $\gamma$):
$$
A_\nu=\sup_{g\in B^+_{L^{\gamma'}}} \int_\T\int_{\Sigma_w} \frac{1}{(1-|z|^2)^\eta}\, g(w)\;d\nu(z) d\lambda(w)=
\sup_{g\in B^+_{L^{\gamma'}}} \int_\D  \frac{{\cal H}(g)(z)}{(1-|z|^2)^\frac{\eta}{2}} \;d\nu(z)
$$
where $B^+_{L^{\gamma'}}$ stands for the positive part of the unit ball of $ L^{\gamma'}(\T)$ and
$$
{\cal H}(g)(z)=\frac{1}{(1-|z|^2)^\frac{\eta}{2}}\int_{\{w|\, z\in\Sigma_w\}}\, g(w)\; d\lambda(w)
$$
In the same way,
$$
B_\nu=\sup_{g\in B^+_{L^{\gamma'}}} \int_\T\int_\D \frac{1}{|1-\overline{w}z|^\eta}\, g(w)\;d\nu(z) d\lambda(w)=
\sup_{g\in B^+_{L^{\gamma'}}} \int_\D  \frac{{\cal K}(g)(z)}{(1-|z|^2)^\frac{\eta}{2}} \;d\nu(z)
$$
where
$$
{\cal K}(g)(z)=\int_\T\,\Bigg(\frac{1-|z|^2}{|1-\overline{w}z|^2}\Bigg)^\frac{\eta}{2}\, g(w)\; d\lambda(w).
$$

{\bf Claim:} there exists some $C>0$ such that, for every positive function 
$g$ and every $z$ in $\D$, $${\cal K}(g)(z)\le C {\cal H}\big({\cal M}(g)\big)(z)$$
where $\dis{\cal M}(g)$ is the Hardy-Littlewood maximal function associated to $g$.

We postpone the proof of the claim and we now get (via Fubini):
$$
\int_\D  \frac{{\cal K}(g)(z)}{(1-|z|^2)^\frac{\eta}{2}} \;d\nu(z)\lesssim 
\int_\D  \frac{{\cal H}\big({\cal M}(g)\big)(z)}{(1-|z|^2)^\frac{\eta}{2}} \;d\nu(z)
\lesssim\int_\T{\cal M} (g)(z)\int_{\Sigma_w} \frac{1}{(1-|z|^2)^\eta}\;d\nu(z)\;d\lambda(w)\,.
$$
Thanks to H\"older inequality, we get 
$$
\int_\D  \frac{{\cal K}(g)(z)}{(1-|z|^2)^\frac{\eta}{2}} \;d\nu(z)\lesssim 
A_\nu\|{\cal M} (g)\|_{\gamma'}\lesssim A_\nu\|g\|_{\gamma'}\,,
$$
since $\gamma'>1$. Passing to the supremum over $g$, we get
$$
B_\nu\lesssim A_\nu
$$
up to constants depending on $\gamma$ and $\eta$ only.
\medskip

Now we prove the Claim. 
We write $z=r\e^{i\theta}$ with $\theta\in\R$ and $0\le r<1$. 
Point out that $z$ belongs to $\Sigma_w$ with $w=\e^{it}$, as soon as 
$t\in\big(\theta-c(1-r),\theta+c(1-r)\big)$, where $c$ is a (numerical) constant. Therefore, identifying the intervals $I$ with the arcs $\{\e^{it}; t\in I\}$, we have
\begin{equation}\label{convol} 
{\cal H}(g)(z)\ge
\frac{1}{(1-|z|^2)^\frac{\eta}{2}}\int_{\mathbb T}\ind_{(\theta-\eps,\theta+\eps)}(w)\,g(w)\;d\lambda(w)
\gtrsim\frac{1}{\eps^\frac{\eta}{2}}\ind_{(-\eps,\eps)}\ast g\big(\e^{i\theta}\big)
\end{equation}
where $\eps=c(1-r)\approx1-|z|^2$.

Besides the  Hardy-Littlewood  maximal function defined by
$$
{\cal M} (g)\big(\e^{i\theta}\big)=\sup_{\e^{i\theta}\in I\subset\T}\frac{1}{\lambda(I)}\int_Ig\,d\lambda, 
$$
we consider also 
$$
{\cal M}_\eps(g)\big(\e^{i\theta}\big)=
\hskip-10pt\biindice{\sup}{\e^{i\theta}\in I}{|I|\ge2\eps}\frac{1}{\lambda(I)}\int_Ig\,d\lambda.
$$ 
It is easy to see that, using (\ref{convol}), that
$$
M:={\cal M}_\eps (g)\big(\e^{i\theta}\big)\le 
{\cal M}(g)\ast\frac{1}{\eps}\ind_{(-\eps,\eps)}\big(\e^{i\theta}\big)=
\eps^{\frac{\eta}{2}-1}\dis{\cal M}(g)\ast\Big(\frac{1}{\eps^\frac{\eta}{2}}\ind_{(-\eps,\eps)}\Big)\big(\e^{i\theta}\big)\lesssim\eps^{\frac{\eta}{2}-1}{\cal H}\big({\cal M}(g)\big)(z) .
$$
Indeed: by definition, there is an interval $I$ with length larger than $2\eps$ 
realizing (almost) the upper bound for the definition of $M$. This interval $I$ contains either  the first half
$(\theta-\eps,\theta)$ or the second half $(\theta,\theta+\eps)$, and any $s$ in the contained half
realizes now ${\cal M}_\eps (g)(\e^{i\theta})\le\dis{\cal M} (g)(\e^{is})$. This is a fortiori true for its mean.
\smallskip

We wish to bound
$$
{\cal K}(g)(z)=P_r^\frac{\eta}{2}\ast g\big(\e^{i\theta}\big)=\int_{(-\pi,\pi)} P_r^\frac{\eta}{2}(t)\, dv(t)
$$ 
where $v$ is the measure  
$$
dv(t)=g\big(\e^{i(\theta-t)}\big)\frac{dt}{2\pi}.
$$
We compute 
$$
{\cal K}(g)(z)=\int_0^{\|P_r\|^\frac{\eta}{2}_\infty} v(\{P_r^\frac{\eta}{2}>x\})\,dx\,.
$$ 
but we already know that $v(I)\le \dis M\lambda(I)$ for every interval containing $0$ with length larger 
than $2\eps$ (by definition). The set $\{P_r>x\}$ is actually an interval $I=(-a,a)\subset(-\pi,\pi)$ 
symmetric with respect to $0$ (thanks to the usual properties of the Poisson kernel:  parity and monotony).
We shall use either the fact $v(I)\le M \lambda(I)$ when $a\ge\eps$, or $v(I)\le2 \eps M $ when $a<\eps$. 
We get 
$$
v(\{P_r^\frac{\eta}{2}>x\})\le M\big(2\eps+\lambda(\{P_r^\frac{\eta}{2}>x\}) \big).
$$

Since ${\eta}>1$, we have 
$$
\int_0^{+\infty}\lambda(\{P_r^\frac{\eta}{2}>x\})\,dx=\big\|P_r \big\|_ \frac{\eta}{2}^\frac{\eta}{2}
\approx \eps^{1-\frac{\eta}{2}}
$$
and we obtain
$$
{\cal K}(g)(z)\le M\Big(2\eps\|P_r\|^\frac{\eta}{2}_\infty+\|P_r\|^\frac{\eta}{2}_\frac{\eta}{2}\Big)
\lesssim M\eps^{1-\frac{\eta}{2}}\lesssim {\cal H}\big({\cal M}(g)\big)(z)
$$
which is the conclusion of the claim.\hfill\cqfd

The next theorem is the main result in this section and it finishes the proof of \eqref{p<2} and completes our characterizations of $r$-summing Carleson embeddings.

\begin{theorem}\label{Mainp<2}
Let  $1<p\le2$ and $\mu$ be a finite measure on the unit disk.
Let $J_\mu\colon H^p\to L^p(\mu)$ be the Carleson embedding, and define
$$
\Phi(\xi)=\int_{\Sigma_\xi}  \frac{1}{\big(1-|z|^2\big)^{1+\frac{p}{2}}}\;d\mu(z)\,,\qquad \xi\in \T\,,
$$
where $\Sigma_\xi$ is the Stolz domain at point $\xi\in\T$. 
Then, for every $r\ge 1$, we have
$$
\pi_r(J_\mu)\approx \pi_2(J_\mu)\approx \|\Phi\|_{L^{2/p}}^{1/p}.
$$
\end{theorem}
 
Let us mention that the case $p=2$ in the previous theorem is already known: the Hilbert-Schmidt norm of a Carleson embedding is equivalent to
$$
\bigg(\int_{\mathbb D}\frac{1}{\big(1-|z|^2\big)}\,d\mu\bigg)^{1/2}\,,
$$
and, using \eqref{hkz}, this is equivalent to
$$
\bigg(\int_\T \int_{\Sigma_\xi}  \frac{1}{\big(1-|z|^2\big)^{2}}\;d\mu(z)\;d\lambda(\xi)\bigg)^{1/2}=
\|\Phi\|_1^{1/2}.
$$

Theorem~\ref{Mainp<2}
will follow from the next more general statement.

\begin{theorem}\label{sommantp,q<2}
Let $1<p\le2$, $1\le q\le2$ and $\mu$ be a finite measure on the unit disk $\D$.
The following assertions are equivalent:

\begin{itemize}
\item The natural injection $J_\mu\colon H^p \to L^q(\mu)$ is $2$-summing

\item The natural injection $J_\mu\colon H^p \to L^q(\mu)$ is $r$-summing for every $r\ge1$.

\item $\dis\xi\longmapsto\int_{\Sigma_\xi}  \frac{1}{\big(1-|z|^2\big)^{1+\frac{q}{2}}}\;d\mu(z)$ belongs to $L^\gamma(\T,d\lambda)$, where $\dis \gamma=\frac{2p}{2p-2q+pq}$ and $\Sigma_\xi$ is the Stolz domain at point $\xi\in\T$.
\end{itemize}

Moreover 
$$\dis \pi_2\bigl(J_\mu\colon H^p\to L^q(\mu)\bigr)\approx \Bigg(\int_\T\Bigg(\int_{\Sigma_\xi}  \frac{1}{\big(1-|z|^2\big)^{1+\frac{q}{2}}}\;d\mu(z)\Bigg)^{\gamma}\;d\lambda(\xi)\Bigg)^{1/q\gamma}.$$
\end{theorem}
\medskip

The heart of the proof of Theorem~\ref{sommantp,q<2} is actually the following proposition which deals with the particular case $q=2$.

\begin{proposition}\label{sommantq=2}
Let $\nu$ be a positive (finite) measure on the unit disk and $1<p\le2$.
TFAE

\begin{enumerate}
\item $\dis J_\nu: H^p\longmapsto L^2(\D,\nu)$ is 2-summing.


\item $\dis\int_\T \biggl( \int_\D \frac{1}{|1-\overline{w}z|^2}\;d\nu(z)  \, \biggr)^{p'/2} \, d\lambda(w)$ is finite. 

\item The Poisson transform, viewed from  $L^p(\T,d\lambda)\longmapsto L^2(\D,\nu)$, is $1$-summing.

\item  $\dis\int_\T \biggl( \int_\D \frac{(1-|z|^2)^2}{ |1-\overline{w}z|^4}\;d\nu(z)  \, \biggr)^{p'/2} \, d\lambda(w)$ is finite. 
\end{enumerate}

Moreover 
 $$\dis \pi_2\bigl(J_\nu\colon H^p\to L^2(\nu)\bigr)\approx \Bigg(\int_\T \biggl( \int_\D \frac{1}{ |1-\overline{w}z|^2}\;d\nu(z)  \, \biggr)^{p'/2} \, d\lambda(w)\Bigg)^{1/p'}.$$
\end{proposition}

{\bf Proof of Proposition~\ref{sommantq=2}. }

$(1)\Rightarrow (2)$ Fix for a while $0<r<1$, define for every $w\in\overline\D$ the function $K_w\in H^p$ by
$$
K_w(z) = \frac{1}{ 1 -r\overline w z}\,, \qquad z\in \overline\D.
 $$ 

Now take $x^*$ in the unit ball of the dual of $H^p$. Considering $H^p$ as a subspace of $L^p(\T)$, there exists $g$ in the unit ball of $L^{p'}(\T)$ such that
$$\langle K_w ,x^*\rangle = \int_\T K_w(z) \overline {g(z)} \,d\lambda(z)=\sum_{n=0}^\infty \overline{\hat g(n)}r^n\overline{w}^n=\overline {g_1(rw)} $$
where $g_1$ is the Riesz projection of $g\in L^{p'}(\T)$ onto $H^{p'}$.


We have then, if $\beta_p$ is the norm of this Riesz projection on $L^{p'}(\T)$, for every $x^*$ in the unit ball 
of the dual of $H^p$,
$$
\Bigl(\int_\T |\langle K_w ,x^*\rangle|^{p'}\,d\lambda(w)\Bigr)^{1/{p'}}
=\Bigl(\int_\T |g_1(rw)|^{p'}\,d\lambda(w)\Bigl)^{1/{p'}}\le \|g_1\|_{H^{p'}}\le \beta_p\|g\|_{p'}$$
hence

$$\dis\sup_{x^\ast\in B_{(H^p)^\ast}}\Bigl(\int_\T |\langle K_w ,x^*\rangle|^{p'}\,d\lambda(w)\Bigr)^{1/{p'}}\le \beta_p\,.$$

Therefore, since $J_\nu$ is ${p'}$-summing (${p'}\ge 2$ and $J_\nu$ is $2$-summing), we have, 
$$\beta_p \pi_2(J_\nu)\ge\beta_p \pi_{p'}(J_\nu) \ge \Bigl(\dis\int_\T \bigl\|K_w\bigr\|_{L^2(\nu)}^{p'}\,d\lambda(w) \Bigr)^{1/{p'}}.$$
So, for every $0<r<1$,
$$\beta_p \pi_2(J_\nu)\ge\Bigg(\int_\T \biggl( \int_\D \frac{1}{ |1-r\overline{w}z|^2}\;d\nu(z)  \, \biggr)^{p'/2} d\lambda(w)\Bigg)^{1/p'}.$$

We get our implication of the statement taking limits when $r\to 1^-$ and using Fatou's Lemma.
\medskip

$(2)\Rightarrow (4)$ is clear as well since $4\big|1-\overline{w}z\big|^2\ge4\big(1-|z|\big)^2\ge \big(1-|z|^2\big)^2$ for any $z\in\D$ and $w\in\T$.
\medskip

$(4)\Rightarrow (3)$. We are interested in the map $f\in L^p(\T,d\lambda)\longrightarrow{\cal P}(f)\in L^2(\D,\nu)$ with 
$$\dis{\cal P}(f)(z)=\int_\T \frac{1-|z|^2}{ |1-\overline{w}z|^2}\,f(w)\; d\lambda(w).$$
{\sl A priori}, it is not obvious that this map is even defined (and bounded). Actually, it appears as the adjoint of the map
$$g\in L^2(\D,\nu)\longrightarrow{\cal Q}(g)\in  L^{p'}(\T,d\lambda)\quad\hbox{with}\quad\dis{\cal Q}(g)(w)=\int_\D \frac{1-|z|^2}{ |1-\overline{w}z|^2}\,g(z)\; d\nu(z).$$
This latter map ${\cal Q}$ is clearly bounded since it is even order bounded: for almost every $w\in\T$
$$\dis\sup_{\|g\|_{L^2(\D,\nu)}\le1}\Bigg|\int_\D \frac{1-|z|^2}{ |1-\overline{w}z|^2}\, g(z)\; d\nu(z)\Bigg|= \biggl( \int_\D \frac{(1-|z|^2)^2}{ |1-\overline{w}z|^4}\;d\nu(z)  \, \biggr)^{1/2} $$
is finite and, as a function of $w\in\T$, it belongs to $L^{p'}(\T,d\lambda)$: it is our hypothesis.

We get that our map ${\cal P}={\cal Q}^\ast$ is defined and bounded as well (with same norm). Moreover Lemma~\ref{LemOB} (point out that ${\cal P}^\ast={\cal Q}$) implies that ${\cal P}$ is $1$-summing.

\medskip

$(3)\Rightarrow (1)$ is clear  by restriction. Indeed, for every $f\in H^p$, we have ${\cal P}(f)=f$.
\hfill\cqfd
\bigskip

\medskip

For the proofs we need also the following lemmas:

\begin{lemma}\label{lemINF}
Let $\sigma>0$, $(\Omega,\Sigma,\mu)$ be a measure space and $h\colon \Omega\to [0,+\infty)$ be a measurable function. Then  
$$
\inf\biggl\{\int_\Omega \frac{h}{ F}\, d\mu : F\in L^\sigma(\mu), F\ge 0, \int_\Omega F^\sigma\,d\mu\le 1\;\biggr\}= 
\Bigr(\int_\Omega h^{\sigma/(\sigma+1)}\, d\mu\Bigl)^{(\sigma+1)/\sigma}.
$$
\end{lemma}

{\bf Proof of Lemma~\ref{lemINF}.} Take $F\ge 0$ in the unit ball of $L^\sigma(\mu)$. Observing that 
$(\sigma+1)/\sigma$ and $\sigma+1$ are conjugate exponents, we have by H\"older,
$$
\int_\Omega h^{\sigma/(\sigma+1)}\,d\mu \le \Bigl(\int_\Omega (h/F)\, d\mu\Bigr)^{\sigma/(\sigma+1)}
\Bigl(\int_\Omega F^\sigma\, d\mu\Bigr)^{1/(\sigma+1)} \le \Bigl(\int_\Omega (h/F)\, d\mu\Bigr)^{\sigma/(\sigma+1)}\,.
$$
In consequence we have, for every $F\ge 0$ in the unit ball of $L^\sigma(\Omega,\mu)$,
$$\Bigr(\int_\Omega h^{\sigma/(\sigma+1)}\, d\mu\Bigl)^{(\sigma+1)/\sigma}\le \int_\Omega h/F \, d\mu$$
 Taking infimum we get one inequality.
\bigskip

To prove the other inequality, we may (and do) assume that $\|h\|_{L^{\sigma/(\sigma+1)}(\mu)}$ is finite. Now taking
$F_0=\lambda h^{1/(\sigma+1)}$, for $\lambda >0$, we see that $F_0\in L^\sigma(\mu)$. 
If we put $\beta=\|h\|_{L^{\sigma/(\sigma+1)}(\mu)}$, 
let us adjust $\lambda$ to get $\|F_0\|_\sigma=1$. We should have
$$
1= \lambda^\sigma \int_\Omega h^{\sigma/(\sigma+1)}\, d\mu = \lambda^\sigma\beta^{\sigma/(\sigma+1)}\,,
$$
and therefore $\lambda=\beta^{-1/(\sigma+1)}$. With this choice of $F_0$, we have
$$
\int_\Omega h/F_0\, d\mu = \lambda^{-1}\int_\Omega h/h^{1/(\sigma+1)} \,d\mu = 
\beta^{1/(\sigma+1)} \int_\Omega h^{\sigma/(\sigma+1)} \,d\mu =
\beta^{1/(\sigma+1)} \beta^{\sigma/(\sigma+1)} = \beta\,.
$$
Then the infimum in the statement is less or equal than 
$$\dis\beta=\Bigr(\int_\Omega h^{\sigma/(\sigma+1)}\, d\mu\Bigl)^{(\sigma+1)/\sigma}$$ and the lemma follows.\cqfd

The following result is probably well known from the specialist, nevertheless we have no explicit reference. We state it and prove it for the convenience of the reader. In the statement, we take the convention that $\dis\frac{0}{0}=0$. 

\begin{lemma}\label{maurey}
Let $1\le q< 2$ and let $s>1$ be such that $1/s + 1/2 = 1/q$. Let $X$ be a Banach space, and
$T\colon X \to L^q(\mu)$ a bounded operator.
The necessary and sufficient condition for $T$ to be a $2$-summing operator is that
there exists $F\in L^s(\mu)$, with $F\ge0$ $\mu$-a.e., such that  
$T\colon X\to L^2(\nu)$ is well defined and $2$-summing, where the measure $\nu$ is the measure defined by
$$
d\nu(z)= \frac{1}{ F(z)^2}\,d\mu(z)\,.
$$
Moreover, we have
$$
\pi_2\bigl(T\colon X\to L^p(\mu)\bigr) \approx \inf\Bigl\{\pi_2\bigl(T\colon X\to L^2(\nu)\bigr)
: d\nu=d\mu/F^2, F\ge0, \int F^s\,d\mu \le 1
\Bigr\}.
$$
\end{lemma}

{\bf Proof of Lemma~\ref{maurey}.} 
Suppose first that $F\in L^s(\mu)$, $F\ge 0$, and that $T\colon X\to L^2(\nu)$ is well defined and $2$-summing, for $d\nu=d\mu/F^2$. Now, we claim that $L^2(\nu) \subset L^p(\mu)$, which is clear since $g\mapsto F\cdot g$ defines actually a multiplier from  $L^2(\mu)$ to $L^p(\mu)$. Writing this for $g=h/F$, the claim is proved. This yields 
$$\pi_2\bigl(T\colon X\to L^p(\mu) \bigr)\le \|F\|_{L^s(\mu)} \pi_2\bigl(T\colon X\to L^2(\nu) \bigr)\,.$$
\medskip

For the converse implication we will use Maurey's Factorization Theorem. 
Suppose that $T\colon X\to L^p(\mu)$ is $2$-summing. Then by Pietsch's factorization, 
there exists a Hilbert space $H$, and
two operators 
$$
S\colon X\to H 
\qquad \hbox{ and } \qquad  R\colon H\to L^p(\mu)\,,
$$
such that $S$ is $2$-summing and $\pi_2(T)=\pi_2(S)$, $\|R\|\le 1$, and $T=R\circ S$.
\bigskip

Given any (finite) family $\{h_i : i\in I\}$ in the unit ball of $H$, 
and any family $\{\alpha_i :i\in I\}$ of real numbers
we have, for $\{r_i : i\in I\}$ a Rademacher family defined on $(\Omega,\Prob)$,
$$
\Bigl( \int \Bigl(\sum_i |\alpha_i Rh_i|^{2}\Bigr)^{p/2} \,d\mu\Bigr)^{1/p} =
\Bigl( \int \Bigl( \int_\Omega \Bigl|\sum_i \alpha_i r_i(\omega) Rh_i\Bigr|^{2}
\,d\Prob(\omega)\Bigr)^{p/2} \,d\mu\Bigr)^{1/p} $$

which is less than
$$c_p    \Bigl( \int \int_\Omega \Bigl|\sum_i \alpha_i r_i(\omega) Rh_i\Bigr|^{p}
\,d\Prob(\omega) \,d\mu\Bigr)^{1/p}$$
where $c_p$ denotes the constant in Khintchin's inequality.

Now using Fubini's theorem and the boundedness of $R$, we get
$$\Bigl( \int \Bigl(\sum_i |\alpha_i Rh_i|^{2}\Bigr)^{p/2} \,d\mu\Bigr)^{1/p}\le
c_p\Bigl( \int_\Omega \Bigl\|\sum_i \alpha_i r_i(\omega) h_i\Bigr\|^{p}_{H}
\,d\Prob(\omega)\Bigr)^{1/p}$$
and since $p\le2$, we have 
$$\Bigl( \int \Bigl(\sum_i |\alpha_i Rh_i|^{2}\Bigr)^{p/2} \,d\mu\Bigr)^{1/p}\le c_p \Bigl( \int_\Omega \Bigl\|\sum_i \alpha_i r_i(\omega) h_i\Bigr\|^{2}_{H}
\,d\Prob(\omega)\Bigr)^{1/2} \le  c_p \Bigl(\sum_i |\alpha_i|^2\Bigr)^{1/2} .$$

We can therefore apply Th\'eor\`eme 2 in the page 12 of Maurey's book \cite{Ma} to the subset  
$$
\Bigl\{ \frac{1}{ c_p}Rh : h\in B_H\Bigr\}
$$
of $L^p(\mu)$, to get a function $F_0\ge 0$, such that $\int F_0^s\, d\mu\le 1$, and
$$
\Bigl(\int |Rh/F_0|^2 \, d\mu \Bigr)^{1/2}\le c_p\,, \qquad\hbox{for all $h\in B_H$.}
$$
Namely, if $d\nu=d\mu/F_0^2$, we have proved that $\bigl\|R\colon H\to L^2(\nu)\bigr\|\le c_p$, 
and consequently,
since $T=R\circ S$, we have $T\colon X\to L^2(\nu)$ is well defined, $2$-summing and
$$
\pi_2\bigl(T\colon X\to L^2(\nu)\bigr) 
\le c_p\pi_2(S)=c_p \pi_2\bigl(T\colon X\to L^p(\mu)\bigr)\,.
$$\hfill\cqfd
\medskip

{\bf Proof of Theorem~\ref{sommantp,q<2}.}
The two first assertions are equivalent since the spaces $H^p$ and $L^q$ have cotype 2 (see \cite{DJT}, cor. 3.16). 
\medskip

Let us treat the particular case $q=2$. This is mainly contained in Prop.~\ref{sommantq=2} although the conclusion of Th.~\ref{sommantp,q<2} involves an integral (over some Stolz domain) of different nature compared to the integrals in Prop.~\ref{sommantq=2}. The equivalence of these integrals is the conclusion of Prop.~\ref{equivInt} in the case $\gamma=1$ and $\eta=2$.
\medskip

Now, we focus on the case $q<2$. Applying Lemma~\ref{maurey} and Prop.~\ref{sommantq=2}, we know that 

$$
\pi_2^2\bigl(J_\mu\colon H^p\to L^q(\mu)\bigr) \approx \inf\Bigg\{\dis\int_\T \biggl( \int_\D \frac{1 }{ |1-\overline{w}z|^2}\,\frac{d\mu(z)}{F^2(z)}  \, \biggr)^{p'/2} \, d\lambda(w)\,\Big|\; F\in B^+_{L^s(\D,\mu)}\Bigg\}^{\frac{2}{p'}}
$$
where $\dis\frac{1}{s}+\frac{1}{2}=\frac{1}{q}$ and $\dis B^+_{L^s(\D,\mu)}$ stands for the positive part of the unit ball of $L^s(\D,\mu)$. 
\smallskip

Replacing $F^2$ by $f$ (and then $s$ by $\sigma=s/2$) and linearizing, we get
$$
\pi_2^2\bigl(J_\mu\colon H^p\to L^q(\mu)\bigr) \approx \dis\inf_{f\in B^+_{L^\sigma(\D,\mu)}}\Bigg(\sup_{g\in B^+_{L^t(\T,d\lambda)}} \int_\T  \int_\D \frac{g(w)}{ |1-\overline{w}z|^2}\,\frac{d\mu(z)}{f(z)} \, d\lambda(w)\Bigg)$$
where $\dis\sigma=\frac{s}{2}=\frac{q}{2-q}$ and $\dis t=\Big(\frac{p'}{2}\Big)'=\frac{p}{2-p}\cdot$
\medskip

{\sl Claim.} We claim now that
$$A:=\inf_{f\in B^+_{L^\sigma(\D,\mu)}}\Bigg(\sup_{g\in B^+_{L^t(\T,d\lambda)}} \int_\T  \int_\D \frac{g(w)}{ |1-\overline{w}z|^2}\,\frac{d\mu(z)}{f(z)} \, d\lambda(w)\Bigg)$$
is equal to $$B:=\sup_{g\in B^+_{L^t(\T,d\lambda)}}\Bigg( \inf_{f\in B^+_{L^\sigma(\D,\mu)}}\int_\T  \int_\D \frac{g(w)}{ |1-\overline{w}z|^2}\,\frac{d\mu(z)}{f(z)} \, d\lambda(w)\Bigg)$$ 

Indeed, $A\ge B$ is obvious and the other inequality is a consequence of the convexity underlying the quantities above. More precisely, we apply the Ky Fan's lemma to the family of functionals $$M=\big\{\Phi_g\,|\; g\in B^+_{L^t(\T,d\lambda)}\big\}$$ where, for $f\in  B^+_{L^\sigma(\D,\mu)}$,
$$\Phi_g(f)=\int_\T  \int_\D \frac{g(w) }{ |1-\overline{w}z|^2}\,\frac{d\mu(z)}{f(z)} \, d\lambda(w)$$

The set $C=B^+_{L^\sigma(\D,\mu)}$, equipped with weak topology, is convex and compact ($\sigma>1$) and $M$ appears as a set of functions on $C$. Since $C$ is convex and the mappings  $g\mapsto \Phi_g$ are linear, the set $M$ itself is convex. 

We assume that $B$ is finite (else the inequality $A\le B$ is trivial) and we are going to see that the three conditions of the Ky Fan's lemma (as stated in \cite{DJT} p.190) are verified. Namely:
\begin{enumerate}[(a)]
\item Each $\Phi\in M$ is convex and lower semi-continuous.

\item If $\Psi\in conv(M)$, there exists some $\Phi\in M$ with $\Psi(x)\le \Phi(x)$ for all $x\in C$.

\item There is an $r\in\R$ such that each $\Phi\in M$ has a value less than $r$.
\end{enumerate}

Fix $\eps>0$ and consider $r=B+\eps$, so that (c) is verified by definition of $B$. 

Conditions (b) is obviously verified since $M$ is convex. 

The functionals $\Phi_g$ are convex (thanks to the convexity of $x\in(0,+\infty)\mapsto 1/x$) so the first part of (a) is verified and we only have to check now that they are also lower semi-continuous. Fixing $\lambda>0$ and $g\in B^+_{L^t(\T,d\lambda)}$, we wish to prove that the convex set $K=\{f\in C\,|\; \Phi_g(f)\le\lambda\}$ is closed in the weak topology. It suffices to prove that $K$ is closed for the strong topology. Take a sequence $(f_n)$ in $K$, converging to some $f$. Up to an extraction, we can assume that $f_n$ is also pointwise converging a.e. to $f$. Hence, by Fatou's lemma,
$$\Phi_g(f)\le\underline{\lim}\Phi_g(f_n)\le\lambda.$$

The conclusion of Ky Fan's lemma says that there exists some $f_0\in C$ such that $$\dis\sup_{g\in B^+_{L^t(\T,d\lambda)}}\Phi_g(f_0)\le r=B+\eps .$$ A fortiori, $A\le B+\eps$ and the claim is proved.
\medskip

{\sl End of the proof.} We get 
$$\pi_2^2\bigl(J_\mu\colon H^p\to L^q(\mu)\bigr) \approx \dis \sup_{g\in B^+_{L^t(\T,d\lambda)}}\Bigg( \inf_{f\in B^+_{L^\sigma(\D,\mu)}}\int_\T  \int_\D \frac{g(w)}{|1-\overline{w}z|^2}\,\frac{d\mu(z)}{f(z)} \, d\lambda(w)\Bigg)$$
which can be written (thanks to Fubini's theorem)
$$\pi_2^2\bigl(J_\mu\colon H^p\to L^q(\mu)\bigr) \approx \dis \sup_{g\in B^+_{L^t(\T,d\lambda)}}\Bigg( \inf_{f\in B^+_{L^\sigma(\D,\mu)}} \int_\D \frac{\mathscr{P}(g)(z)}{ 1-|z|^2}\,\frac{d\mu(z)}{f(z)}\Bigg)$$
where $\mathscr{P}(g)$ is the Poisson transform of $g$.

Applying Lemma~\ref{lemINF} and replacing $\sigma$ by its value, we obtain 
$$\pi_2^2\bigl(J_\mu\colon H^p\to L^q(\mu)\bigr) \approx \dis \sup_{g\in B^+_{L^t(\T,d\lambda)}} \Bigg\|\frac{\mathscr{P}(g) }{ 1-|z|^2} \Bigg\|_{L^{\frac{q}{2}}(\D,d\mu)}$$
but it means that the Poisson transform 
$$\dis g\in L^t(\T,d\lambda)\longmapsto \mathscr{P}(g)\in L^{\frac{q}{2}}\Big(\D,\frac{d\mu}{(1-|z|^2)^{\frac{q}{2}}}\Big)$$
is bounded, with a norm equivalent to $\pi_2^2\bigl(J_\mu\colon H^p\to L^q(\mu)\bigr)$.

Since $t>1$ (because $p>1$), the boundedness of this Poisson transform is equivalent to the fact that the Hardy space $H^t$ is sent into  $\dis L^{\frac{q}{2}}\Big(\D,\frac{d\mu}{(1-|z|^2)^{\frac{q}{2}}}\Big)$.

Let us check that we are actually working with a finite measure $\dis\frac{d\mu}{(1-|z|^2)^{\frac{q}{2}}}\cdot$ Indeed, if we assume that the embedding is absolutely summing, then this measure is finite thanks to (\ref{ncpq}) (cf Prop.~\ref{propnec}). Conversely, if we assume that the integral condition is fulfilled: $\dis\xi\longmapsto\int_{\Sigma_\xi}  \frac{1}{\big(1-|z|^2\big)^{1+\frac{q}{ 2}}}d\mu(z)$ belongs to $L^\gamma(\T,d\lambda)$, in particular it belongs to $L^1(\T,d\lambda)$:
$$\int_{\mathbb T}\int_{\Sigma_\xi}  \frac{1}{\big(1-|z|^2\big)^{1+\frac{q}{2}}}d\mu(z)\;d\lambda<\infty $$
which implies that (via Fubini)
$$\int_{\mathbb D}  \frac{1}{\big(1-|z|^2\big)^{\frac{q}{ 2}}}\;d\mu(z)\lesssim\int_{\mathbb D}\int_{\{\xi\in\mathbb T; z\in\Sigma_\xi\}} \frac{1}{\big(1-|z|^2\big)^{1+\frac{q}{ 2}}}\;d\lambda\;d\mu(z)<\infty$$
and the measure is finite.

We can now apply directly a work due to Blasco-Jarchow (\cite{BJ}), following a former work of Luecking (\cite{Lue}). The reader can check the statement of Th.1.3 in \cite{BJ}, but see mainly the proof of Th.2.2 (\cite{BJ}) where it is actually proved that for any finite positive measure $\nu$, carried by $\D$, we have for $a>b>0$:
$$\big\|H^a\hookrightarrow L^b(\D,\nu)\big\|\approx\Big\|\xi\longmapsto \int_{\Sigma_\xi}\frac{d\nu}{1-|z^2|}\,\Big\|_{L^c({\mathbb T})}^{\frac{c}{b}} $$
where $c=\dis\frac{a}{a-b}$ and the underlying constants depend on $a$ and $b$.

The conclusion follows.\cqfd

\section{Applications, examples and remarks.}

The first application deals with one of the most famous injection in function theory. It is well known that the (formal) identity from the Hardy space $H^p$ to the Bergman space ${\cal B}^q$ is defined (and bounded) if and only if $q\le2p$, and that it is compact if and only if $q<2p$. It is then a natural question to decide when it is absolutely summing. This will be the aim of Th.~\ref{BergInjSum}. First, we state as a lemma the particular case $p=1$. 

\begin{lemma}\label{BergInjSum_p=1}
Let $q<2$. The injection from the Hardy space $H^1$ to the Bergman space ${\cal B}^q$ is $1$-summing.
\end{lemma}

{\bf Proof.} First we point out that we are in the framework of Carleson's embeddings: our measure $\mu$ is here the area measure ${\cal A}$. We shall use the notations of section 3. In particular, ${\cal A}_N$ is the area measure ${\cal A}$ restricted to the corona ${\cal G}_N$.

Now, we observe that for every fixed $N\ge1$, the $1$-summing norm of the injection ${T}_{1,N}$ is summable. Indeed, for any $p>1$ (but we shall specify our choice of $p$ below), it factorizes through $${H_{1,N}}\stackrel{Id}{\longmapsto} {H_{p,N}}\stackrel{Id}{\longmapsto}{H_{p}}\stackrel{Id}{\longmapsto}L^q(\D,{\cal A}_N).$$
The identity viewed from ${H_{1,N}}$ to ${H_{p,N}}$ has operator norm less than $N^{1/{p'}}$. The $1$-summing norm of the last factor is majorized by $$\dis \Big(\sup_\xi {\cal A}_N(\Sigma_\xi)N^{1+\frac{q}{2}}\Big)^{\frac{1}{q}}$$ thanks to  Th.~\ref{sommantp,q<2} applied to the measure ${\cal A}_N$. Now the rotation invariance of the measure and the geometric properties of the Stolz domain gives that $\dis{\cal A}_N(\Sigma_\xi)={\cal A}_N(\Sigma_1)\le c\frac{1}{N^2}$ for some numerical $c>0$. We get
$$\pi_1\big({T}_{1,N}\big)\lesssim N^{-\eps}$$
where $\dis\eps=\frac{1}{q}-\frac{1}{2}-\frac{1}{p'}$ which turns out to be positive for a suitable choice of $p>1$.

Now we use the same trick than in Th.\ref{CourSomGen} (comparing $(b)$ and $(a)$). Of course, the Riesz transform is not uniformly bounded anymore but for every $f\in H^1$, still using the Riesz projection $\pi_m(f)$ of $f$ on the space $Z_m$, spanned by the $z^k$ when $k$ runs over $\{mN,\ldots,(m+1)N-1\}$, the norm $\dis\|y_m\|_{H^1}$ is now bounded by $\log(N)\|f\|_{H^1}$ (up to a constant). Now the sequel of the argument follows the same lines and we get
$$\dis\pi_1\big(J_{{\cal A}_N}\big)\lesssim\log(N)\pi_1\big({T}_{1,N}\big)\lesssim\log(N)N^{-\eps}.$$

Specifying for values $N=2^n$, we have that $\dis\pi_1\big(J_{{\cal A}_{2^n}}\big)$ is summable. At last, by the triangular inequality $J_{\cal A}: H^1\to L^q(\D,{\cal A})$ is $1$-summing.\cqfd
\medskip

\begin{theorem}\label{BergInjSum}
Let $p,q\ge1$ with $q\le 2p$. The injection from the Hardy space $H^p$ to the Bergman space ${\cal B}^q$ is $r$-summing for some $r\ge1$ if and only if $q<max(2,p)$. 
Moreover
\begin{itemize}
\item When  $q<2$, this operator is $1$-summing.

\item When $2\le q<p$, this operator is $r$-summing exactly for every $r$ such that $\dis\frac{1}{p}+\frac{1}{r}<\frac{2}{q}\cdot$
\end{itemize}\end{theorem}

{\bf Proof.} We shall use that the injection from the Hardy space $H^2$ to the Bergman space ${\cal B}^2$ is not $2$-summing. Indeed, reasoning with Taylor coefficients, it would mean that the diagonal operator from $\ell^2$ to itself with diagonal entries $\frac{1}{\sqrt{n+1}}$ would be Hilbert-Schmidt but the eigenvalues do not belong to $\ell^2$. Another way to prove it is to use the remark just after Prop.\ref{propnec}.

In the case $1\le p\le 2$, if the injection from the Hardy space $H^p$ to the Bergman space ${\cal B}^2$ were $r$-summing for some $r\ge1$, then {\it a fortiori} the injection from the Hardy space $H^2$ to the Bergman space ${\cal B}^2$ would be $r$-summing hence $2$-summing, which was just shown to be false. Hence the injection from the Hardy space $H^p$ to the Bergman space ${\cal B}^2$ is not $r$-summing for any $r\ge1$. 

For every $q<2$, the injection from $H^1$ to the Bergman space ${\cal B}^q$ is $1$-summing, by Lemma \ref{BergInjSum_p=1}. A fortiori, the injection from $H^p$ to the Bergman space ${\cal B}^q$ is $1$-summing for any $p\ge1$ (just factorize through the injection from $H^1$ to ${\cal B}^q$). This settles the case $q<2$.

Now, we assume that $q\ge2$. When $q\ge p$, the $r$-summingness would imply the $\gamma$-summingness for $q=p$ and $\gamma=\max(p,r)$. Then our characterization implies that it should be order bounded (from $H^p$ to $L^p(\D,{\cal A})$) which is not. So now on, we assume that $q<p$ and it is easy  to see that the injection  $H^p$ to the Bergman space ${\cal B}^q$ is $q$-summing since it is order bounded.  

Let us focus on the case $2\le q<p$ and assume first that the injection is $r$-summing  for some $r> p'$. Then we use Lemma~\ref{LemAp'<r<2} and the computation made after. The function $\Phi$ satisfies that $$\dis\sup_{\|g\|_{(H^p)^\ast}\le1}\int_0^{+\infty}|\langle \Phi(t),g\rangle|^r\,dt\le C $$
for some $C>0$. Hence the $r$-summingness of the injection from $H^p$ to ${\cal B}^q$ then implies that $\dis\int_0^{+\infty}\|\Phi(t)\|_{{\cal B}^q}^r\,dt$ is bounded. But reproducing the computations in our framework, we get:
$$ \dis\int_0^{+\infty}\|\Phi(t)\|^r_{L^q({\cal A})}\,dt\ge\sum_{n\ge1}\sum_{j=0}^{2^{n}-1}2^{-n}\frac{2^{rn(\frac{1}{p}+\frac{1}{r})}}{2^r}{\cal A}\big(R_{n,j}\big)^{r/q}\approx\sum_{n\ge1}\frac{2^{rn(\frac{1}{p}+\frac{1}{r})}}{2^r}4^{-nr/q}.$$ 
So we necessarily have $\dis\frac{1}{p}+\frac{1}{r}<\frac{2}{q}\cdot$ We have to show that the situation $r\le p'$ is not possible. But if our operator were $r$-summing for some $r\le p'$ then it would be $s$-summing for any $s> p'$ and our previous condition applies: $\dis\frac{1}{p}+\frac{1}{s}<\frac{2}{q}$ which gives (passing to the limit) $\dis q\le2$. 

So it remains to justify that $H^p\hookrightarrow{\cal B}^2$ cannot be $p'$-summing. Indeed, assume the contrary. The formal Riesz projection viewed from $C(\T)$ to $H^{p}$ is $p$-summing since $C(\T)\hookrightarrow L^p(\T)$ is $p$-summing and the Riesz projection is bounded from $L^p(\T)$ to $H^p$ ($p>1$ here). Hence the composition with our injection is now $1$-summing. Using for instance the Pietsch theorem and the translation invariance (on the torus), equivalently the rotational invariance, of this Riesz projection from $ C(\T)$ to ${\cal B}^2$, the Pietsch measure can be chosen as the Haar measure on the torus. We get for every trigonometric polynomial $f$:
$$\dis\Big(\sum_{n\ge0}\frac{|\hat f(n)|^2}{n+1}\Big)^\frac{1}{2}\lesssim\|f\|_{L^1}.$$
Testing now this latter inequality for instance on  the Poisson kernels $P_\rho$, we get that all the sums $\dis\sum_{n\ge0}\frac{\rho^{2n}}{n+1}$ should be bounded independently of $\rho\in(0,1)$, which is false.

Now, for $2\le q<p$, fix $r\le q$ satisfying the condition $\dis\frac{1}{p}+\frac{1}{r}<\frac{2}{q}$ and consider the Poisson kernel $P_\rho$ associated to $\dis\rho=\sqrt{1-2^{-(n+1)}}$. We shall first work with a dyadic corona. The operator $f\in H^p\mapsto f(\rho z)\in H^\infty$ is bounded with norm less that $\dis2^{n/p}$ (up to a constant). The injection from $H^\infty$ to $H^r$ is $r$-summing (with norm $1$). The operator $f\in H^r\mapsto f(\rho z)\in H^q$ is bounded with norm less that $\dis 2^{n(\frac{1}{r}-\frac{1}{q})}$ and at last the operator $g\in  H^q\mapsto g\big(\rho^{-2} z\big)\in L^q(\Gamma_{n},{\cal A})$ is bounded with norm less than $\dis 2^{-n/q}$. By composition, we get that the operator $J_{{\cal A}_{2^n}}$ is $r$-summing with an $r$-summing norm less than $\dis2^{n\theta}$ (up to a constant) where $\dis\theta=\frac{1}{p}+\frac{1}{r}-\frac{2}{q}<0$. Since the series converges, we get that our injection is $r$-summing.\cqfd

\bigskip

The following theorems state that we can separate the different classes of $r$-summing on $H^p$ spaces (as soon as it is possible, according to the values of $r$ or $p$) using only  composition operators:

\begin{theorem}\label{SeparerSelonR}
Fix $p>2$.  For every $r,s\ge1$ with $r>\max(s,p')$ and $s<\min(r,p)$, there exists a symbol $\varphi$ such that  $C_\varphi:H^p\to H^p$  is $r$-summing but not $s$-summing.

\end{theorem}
Recall that when $p\le2$, for every $r,s\ge1$, as for any operator on $H^p$ (which has cotype $2$ here), $C_\varphi$ is $r$-summing  if and only if  $\; C_\varphi$ is $s$-summing.
On the other hand, when $p>2$ and $r,s\le p'$, any  $r$-summing operator on $H^p$  is $s$-summing.
At last, we cannot separate $r$-summing operators with composition operators when $r\ge p$ since they all coincide with order bounded composition operators.
\medskip

On the other hand, we have a monotonicity  relatively to $p$ for any Carleson embedding. Moreover, the class of composition operators is large enough to separate the different classes of summing operators:

\begin{theorem}\label{SeparerSelonP}
Fix $r\ge1$.  The mapping ${p>1}\longmapsto$ the class of $r$-summing Carleson embeddings on $H^p$ is non increasing.

More precisely:
\begin{itemize}
\item For any Carleson measure $\mu$: for every $\dis  p_2\ge p_1>1$, when $J_\mu$ is $r$-summing on $H^{p_2}$, then  $J_\mu$  is $r$-summing  on $H^{p_1}$.

\hskip-10pt Moreover

\item Assume that $\dis1<p_1<p_2$.  If $p_1<2$, we consider any value of $r\ge1$ and if $p_1\ge2$, we consider the values of $r$ such that $p_2>r\ge1$. Then there exists a symbol $\varphi$ such that  $C_\varphi$  is $r$-summing on $H^{p_1}$ but not on $H^{p_2}$.
\end{itemize}
\end{theorem}

Now we prove the two aforementioned theorems. We have to separate the cases since our characterization depends on the values of $(p,r)$. Actually, we could have made a more general (and synthetic)  discussion  which would rely on the distinction of cases according to where stands $(p,r)$ in a diagram drawn in $(1,+\infty)\times(1,+\infty)$ with the curves $p\rightarrow r=p$ and $\dis p\rightarrow r=p'=\dis\frac{p}{p-1}$; with this viewpoint the line with gradient equal to the quotient $\frac{r}{p}$ plays a key role. Nevertheless, we chose to focus on the distinctions fixing either $r$ or $p$.
\medskip

{\bf Proof of Th. \ref{SeparerSelonR}.} Fix $p>2$. It suffices to prove the case $p\ge r>s\ge p'$ since we know that  $C_\varphi$  is $p'$-summing if and only $C_\varphi$  is $1$-summing on one hand. And on the other hand, we know that $C_\varphi$  is $r$-summing if and only $C_\varphi$  is $p$-summing when $r\ge p$.
 
We first choose some $\rho$ such that $s<\rho<r$. So we have $\dis2r/p> 2\rho/p\,\cdot$ Now we consider a symbol $\varphi$ introduced in Th.5.1 \cite{LLQR3} in order to  satisfy $\dis C_\varphi\in{\cal S}_\frac{2r}{p}(H^2)\setminus{\cal S}_\frac{2\rho}{p}(H^2)$. But we already noticed that $C_\varphi$ is an $r$-summing composition operator on $H^p$ if and only if $\dis C_\varphi$ belongs to the Schatten class $\dis{\cal S}_\frac{2r}{p}(H^2)$. Hence $C_\varphi$ is $r$-summing but not $\rho$-summing, a fortiori not $s$-summing.\cqfd
\medskip

{\bf Proof of Th.\ref{SeparerSelonP}.} Fix $r\ge1$ and $p_2>p_1>1$.
\begin{itemize}

\item If $1<p_1<p_2\le2$, we use the argument given in the proof of Th.\ref{sommantp,q<2}. Assume that $\dis J_\mu\colon H^{p_2}\to L^{p_2}(\D,\mu)$  absolutely summing, then
$$\dis \pi_2^2\bigl(J_\mu\bigr)\approx\int_\T\Bigg(\int_{\Sigma_\xi}  \frac{1}{\big(1-|z|^2\big)^{1+\frac{p_2}{2}}}d\mu(z)\Bigg)^{\frac{2}{{p_2}}}\;d\lambda(\xi)\qquad\hbox{is finite}$$
which is equivalent to the fact that $\dis H^\frac{p_2}{2-p_2}$ is sent into $\dis L^\frac{p_2}{2}\Big(\D,\frac{d\mu}{(1-|z|^2)^{\frac{p_2}{2}}}\Big)$. But this is equivalent also to the fact that 
$\dis H^2$ is sent into $\dis L^{2-p_2}\Big(\D,\frac{d\mu}{(1-|z|^2)^{\frac{p_2}{2}}}\Big)$ (see \cite{BJ}), say with norm $K_2$.

Now fix a function $f$ in the unit ball of $H^2$, for every $z\in\D$, the evaluation $|f(z)|$ is majorized by $\dis\frac{1}{(1-|z|^2)^{\frac{1}{2}}}\cdot$ Therefore, using $p_2-p_1>0$, we point out that
$$\int_\D |f(z)|^{2-p_1}\frac{d\mu}{(1-|z|^2)^{\frac{p_1}{2}}}\le\int_\D |f(z)|^{2-p_2}\frac{1}{(1-|z|^2)^{\frac{p_2-p_1}{2}}}\cdot\frac{d\mu}{(1-|z|^2)^{\frac{p_1}{2}}}$$

which is equal to $\dis\int_\D |f(z)|^{2-p_2}\frac{d\mu}{(1-|z|^2)^{\frac{p_2}{2}}}\le K_2$.

Hence $\dis H^2$ is sent into $\dis L^{2-p_1}\Big(\D,\frac{d\mu}{(1-|z|^2)^{\frac{p_1}{2}}}\Big)$ and $\dis J_\mu\colon H^{p_1}\to L^{p_1}(\D,\mu)$ is absolutely summing.
\medskip

Now, we work with the symbol considered in Lemma 3.7, Lemma 4.3 and the proof of  Th. 4.1. \cite{LLQR3} associated to $t^\beta$ where $\beta=\dis\frac{2}{p_2}\in(0,2]$. With this symbol, the size of the Luecking boxes is controlled as follows: for every $n,j$,
$$\dis \lambda_\varphi(R_{n,j})\approx 2^{-n(1+\frac{p_2}{2})}.$$

We can now apply our Proposition \ref{Encadrep<2}:

On one hand, the series $\dis \sum_{n\ge0}\sum_{j=0}^{2^n-1}\Big(2^{n}\lambda_\varphi(R_{n,j})\Big)^{2/p_2}$ diverge so $C_\varphi$ is not absolutely summing on $H^{p_2}$.

On the other hand, the series $\dis \sum_{n\ge0}\Big(\sum_{j=0}^{2^n-1}\Big(2^{n}\lambda_\varphi(R_{n,j})\Big)^{2/p_1}\Big)^{p_1/2}$ converges so $C_\varphi$ is absolutely summing on $H^{p_1}$.

\item If $1<p_1<2<p_2$ then any $r$-summing $J_\mu$ on $H^{p_2}$ is order bounded on $H^{p_2}$ {\sl a fortiori} on $H^{2}$, which in turn, defines an $r$-summing operator on $H^{2}$, hence on $H^{p_1}$, thanks to the previous case. Moreover, the previous case provides us too with an $r$-summing composition operator  on $H^{p_1}$ which is not $r$-summing on $H^{2}$, {\sl a fortiori} it cannot be an $r$-summing operator on $H^{p_2}$.

\item When $2\le p_1<p_2\le r\,$, a Carleson embedding $J_\mu$ is $r$-summing  on $H^{p_2}$ if and only if it is order bounded on $H^{p_2}$, equivalently on $H^{p_1}$, if and only if   $J_\mu$ is $r$-summing  on $H^{p_1}$.

\item If $2\le p_1<p_2$ and $r<p_2$, we split this case in three sub-cases.

\begin{itemize}
\item If $r\ge p_1$, any $r$-summing $J_\mu$ on $H^{p_2}$ is $p_2$-summing on $H^{p_2}$, hence order bounded, which is equivalent to $r$-summing on $H^{p_1}$ by our characterization (\ref{r>p}). On the other hand, there exists a $p_2$-summing composition operator on $H^{p_2}$ ({\sl a fortiori}  $r$-summing  on $H^{p_1}$ by the previous argument), which is not $r$-summing on $H^{p_2}$ by Th.\ref{SeparerSelonR}. 

\item If $p_1'\le r< p_1$, write $s=\dis\frac{rp_2}{p_1}\in(r,p_2)\subset(p_2',p_2)$. Any $r$-summing Carleson embedding  $J_\mu$ on $H^{p_2}$ is $s$-summing. Since we have $\dis\frac{r}{p_1}=\frac{s}{p_2}$, $J_\mu$ is $r$-summing on $H^{p_1}$, by our characterization (\ref{p'<r<p}) (and Lemma~\ref{lemp'<r} if $r=p_1'$). On the other hand, there exists an $s$-summing composition operator on $H^{p_2}$ ({\sl a fortiori} it is  $r$-summing  on $H^{p_1}$ by the previous argument) which is not $r$-summing on $H^{p_2}$ thanks to Th.\ref{SeparerSelonR}. 

\item  If $r< p_1'$, any $r$-summing Carleson embedding $J_\mu$ on $H^{p_2}$ is $p_1'$-summing on $H^{p_2}$, hence $p_1'$-summing on $H^{p_1}$ thanks to the previous case. So $J_\mu$ is actually $1$-summing on $H^{p_1}$ thanks to our characterization (\ref{1<r<p'}), therefore $r$-summing on $H^{p_1}$. On the other hand, thanks to the previous case, there exists a $p_1'$-summing composition operator on $H^{p_1}$ ({\sl a fortiori} it is $r$-summing  on $H^{p_1}$) which is not $p_1'$-summing on $H^{p_2}$. This composition operator cannot be $r$-summing on $H^{p_2}$.

\end{itemize}
\end{itemize}
 




\end{document}